\input amstex

\def\nologo{\let\logo@\empty}

\def\qaq{\quad\hbox{ and }\quad}

\def\qif{\quad\hbox{ if }}

\def\sp{\hbox{ }}

\def\scs{\sp:\sp}

\def\nab{\nabla}

\def\pr{\vskip 5pt \noindent \hbox{\it Proof} \rm \sp \sp}

\def\Q{\Gamma}

\def\Proj{\hbox{\rm Proj}}

\def\dim{\hbox{\rm dim}}

\def\exp{\hbox{\rm exp}}

\def\max{\hbox{\rm max}}

\def\log{\hbox{\rm log}}

\def\PGL{\hbox{\rm PGL}}
\def\Proj{\hbox{\rm Proj}}

\def\Spec{\hbox{\rm Spec}}
\def\Coker{\hbox{\rm Coker}}

\def\Ker{\hbox{\rm Ker}}

\def\Im{\hbox{\rm Im}}
\def\bQ{\Bbb Q}
\def\bZ{\Bbb Z}

\def\bR{\Bbb R}
\def\bP{\Bbb P}

\def\bC{\Bbb C}

\def\2z{\bZ/2\bZ}
\def\2nuz{\bZ/2^\nu\bZ}

\def\cO{{\Cal O}}

\def\cH{{\Cal H}}

\def\Th#1.{\vskip 6pt \medbreak\noindent{\bf Theorem(#1).}}
\def\Pr#1.{\vskip 6pt \medbreak\noindent{\bf Proposition(#1).}}
\def\Lem#1.{\vskip 6pt \medbreak\noindent{\bf Lemma(#1).}}
\def\Rem#1.{\vskip 6pt \medbreak\noindent{\bf Remark(#1).}}
\def\Fact#1.{\vskip 6pt \medbreak\noindent{\bf Fact(#1).}}
\def\Claim#1.{\vskip 6pt \medbreak\noindent{\bf Claim(#1).}}

\def\Conj#1.{\vskip 6pt \medbreak\noindent{\bf Conjecture(#1).}}
\def\Quest#1.{\vskip 6pt \medbreak\noindent{\bf Question(#1).}}
\def\Def#1.{\vskip 6pt \medbreak\noindent{\bf Definition\bf(#1)\rm.}}
\def\Ass#1.{\vskip 6pt \medbreak\noindent{\bf Assumption(#1).}}
\def\Ex#1.{\vskip 6pt \medbreak\noindent{\bf Example(#1).}}
\def\Cor#1.{\vskip 6pt \medbreak\noindent{\bf Corollary(#1).}}

\def\mod{\sp\hbox{\rm mod}\sp}

\def\isom{@>\cong>>}

\def\Xs{X_s}
\def\Zs{Z_s}

\def\scs{\spa:\spa}
\def\Xs{X_s}
\def\Zs{Z_s}

\def\ZZ#1{Z^{[#1]}}
\def\sp{\hbox{ }}
\def\scs{\sp :\sp}

\def\tM{\tilde{M}}
\def\tM'{\tilde{M}'}
\def\otC{\otimes_{\bC}}
\def\ottC{\otimes{\bC}}
\def\otZ{\otimes_{\bZ}}
\def\otZQ{\otimes_{\bZ}{\bQ}}

\def\otO{\otimes_{\cO}}

\def\calX{{\Cal X}}
\def\calZ{{\Cal Z}}

\def\cZS{{\Cal Z}_{S}}


\def\WS#1{\Omega_S^{#1}}
\def\WSd#1{\Omega_{S,d=0}^{#1}}
\def\WXS#1{\Omega_{X/S}^{#1}}
\def\WXSZ#1{\Omega_{X/S}^{#1}(\log Z)}
\def\WXSZgeq#1{\Omega_{X/S}^{\geq #1}(\log Z)}

\def\WXlogZ#1{\Omega_{X}^{#1}(\log Z)}
\def\WX#1{\Omega_{X}^{#1}}

\def\QDZX#1{\bQ_D(#1)_{Z}}
\def\Cone{\hbox{\rm Cone}}
\def\Xan{X_{an}}

\def\NCU{N_{\bC}(U/S)}

\def\MQU{M_{\bQ}(U/S)}
\def\MCU{M_{\bC}(U/S)}


\def\WSd#1{\Omega^{#1}_{S,d=0}}

\def\HO{H_{\cO}}
\def\HQ{H_{\bQ}}

\def\KQ{K_{\bQ}}

\def\JQ{J_{\bQ}}

\def\HQ{H_{\bQ}}

\def\HC{H_{\bC}}

\def\HQX#1#2{\HQ^{#1}(X/S)(#2)}
\def\HCX#1{\HC^{#1}(X/S)}

\def\HCU#1{\HC^{#1}(U/S)}
\def\HCZX#1{\HC^{#1,Z}(X/S)}
\def\HOX#1{\HO^{#1}(X/S)}
\def\HOXpr#1{\HO^{#1}(X/S)_{pr}}

\def\HOU#1{\HO^{#1}(X/S)}

\def\HCU#1{\HC^{#1}(X/S)}

\def\cHDXn#1{\cH^{#1}_D(X/S)(n)}

\def\cHDUn#1{\cH^{#1}_D(U/S)(n)}

\def\cHDZXn#1{{\cH^{#1}_{D,Z}(X/S)}(n)}

\def\QDX#1{\bQ_D(#1)_{X}}
\def\QDU#1{\bQ_D(#1)_{U}}

\def\QDZX#1{\bQ_D(#1)_{(Z,X)}}

\def\HQX#1#2{\HQ^{#1}(X/S)(#2)}
\def\HQXpr#1#2{\HQ^{#1}(X/S)_{pr}(#2)}
\def\HNX#1{H_{\nabla}^{#1}(X/S)}
\def\HNZX#1{H_{\nabla,Z}^{#1}(X/S)}
\def\HQZX#1#2{H^{#1}_{\bQ,Z}(X/S)(#2)}

\def\HOX#1{\HO^{#1}(X/S)}
\def\HOXpr#1{\HO^{#1}(X/S)_{pr}}
\def\HOZX#1{H^{#1}_{\cO,Z}(X/S)}

\def\HOU#1{\HO^{#1}(U/S)}

\def\HHOX#1#2{\HO^{{#1},{#2}}(X/S)}
\def\HHOXpr#1#2{\HO^{{#1},{#2}}(X/S)_{pr}}
\def\HHOZX#1#2{H^{{#1},{#2}}_{\cO,Z}(X/S)}
\def\HHOU#1#2{\HO^{{#1},{#2}}(U/S)}

\def\HCX#1{\HC^{#1}(X/S)}
\def\HCZX#1{H^{#1}_{\bC,Z}(X/S)}

\def\HQZi#1#2{\HQ^{#1}(Z_i/S)(#2)}

\def\HQU#1#2{\HQ^{#1}(U/S)(#2)}

\def\HCU#1{\HC^{#1}(U/S)}



\def\PZr{\Phi^r_{Z/S}}
\def\PZ#1{\Phi^{#1}_{Z/S}}

\def\PUr{\Phi^r_{U/S}}
\def\PU#1{\Phi^{#1}_{U/S}}

\def\FZr{\Psi^r_{Z/S}}

\def\FZ#1{\Psi^{#1}_{Z/S}}

\def\FU#1{\Psi^{#1}_{U/S}}
\def\FXr{\Psi^r_{X/S}}
\def\FXindr{\Psi^r_{X/S,ind}}

\def\CZr{\Lambda^r_{Z/S}}
\def\CXr{\Lambda^r_{X/S}}
\def\CXindr{\Lambda^r_{X/S,ind}}


\def\nabb{\overline{\nab}}

\def\cO{\Cal O}
\def\nabb{\overline{\nab}}

\def\oX#1#2{\rho_X^{#1,#2}}
\def\oS#1#2{\rho_S^{#1,#2}}
\def\oU#1#2{\rho_U^{#1,#2}}
\def\oZ#1#2{\rho_{Z,X}^{#1,#2}}

\def\pZS#1{\phi_{Z/S}^{#1}}

\def\pZSr{\phi_{Z/S}^r}

\def\cZS#1{\psi_{Z/S}^{#1}}
\def\cZSr{\psi_{Z/S}^r}
\def\cUS#1{\psi_{U/S}^{#1}}

\def\cXSr{\psi_{X/S}^r}

\def\tZSr{\tau_{Z/S}^r}
\def\tXSr{\tau_{X/S}^r}

\def\pUSr{\phi_{U/S}^r}
\def\pUS#1{\phi_{U/S}^{#1}}

\def\cUS#1{\psi_{U/S}^{#1}}

\def\cXSr{\psi_{X/S}^r}

\def\tXSr{\tau_{X/S}^r}

\def\San{S_{an}}

\def\otC{\otimes_{\bC}}

\def\nab#1{\nabla_{#1}}
\def\nabb#1{\overline{\nabla}_{#1}}

\def\map#1{@>#1>>}

\def\sE#1#2#3{E^{#1,#2}_{#3}}
\def\sIE#1#2#3{{^{I}E^{#1,#2}_{#3}}}
\def\sIE#1#2#3{{^{I}E^{#1,#2}_{#3}}}
\def\sIIEX#1#2#3{{^{II}E^{#1,#2}_{#3}(X)}}
\def\sIIEZ#1#2#3{{^{II}E^{#1,#2}_{#3}(Z)}}

\def\San{S_{an}}

\documentstyle{amsppt}
\hsize=16cm
\vsize=23cm

\topmatter

\title 
On $K_1$ and $K_2$ of algebraic surfaces
\endtitle

\author{Stefan M\"uller-Stach and Shuji Saito}
\endauthor

\address
{S. M\"uller-Stach: Fachbereich 6, Mathematik und Informatik,
Universit\"at Essen,
45117 ESSEN, Germany}
\endaddress
\address
{e-mail: mueller-stach\@uni-essen.de}
\endaddress
\address
{S. Saito: Graduate School of Mathematics,
Nagoya University, Chikusa-ku, NAGOYA, 464-8602, JAPAN}
\endaddress
\address
{e-mail: sshuji\@msb.biglobe.ne.jp}
\endaddress
\endtopmatter

\document
\NoBlackBoxes

\head \bf Contents 
\endhead

\vskip 10pt\noindent
\roster
\item"\S0"
Introduction
\item"\S1"
Higher Chow groups of normal crossing divisors
\item"\S2"
Normal functions associated to higher cycles on $Z$
\item"\S3"
Proof of the main results
\item"\S4"
Indecomposable parts of infinitesimal invariants
\item"\S5"
An Example in $CH^3(X,2)$
\item"\S6"
An Example in $CH^2(X,1)$
\item"\S7"
Appendix (by Alberto Collino)
\item"\hbox{ }"
References
\endroster

\vskip 10pt

\head \S0. Introduction. \endhead

\vskip 8pt 

Let $X$ be a projective smooth surface over a field $k$ of characteristic zero.
In this paper we study the higher Chow groups $CH^2(X,1)$ and $CH^3(X,2)$ 
using Hodge theoretical methods. They are the most interesting graded pieces
of the Quillen K-groups $K_1(X)$ and $K_2(X)$. 
We recall that these groups are generated by \it higher cycles \rm
that are curves together with sets of rational functions on them 
modulo certain relations arising from tuples of rational functions on $X$. 
More precisely, $CH^2(X,1)$ is the cohomology of the complex
$$ K_2(k(X)) \to \underset{Z\subset X}\to{\bigoplus} k(Z)^* \to 
 \underset{x\in X}\to{\bigoplus} \Bbb Z,$$
where $Z\subset X$ ranges over all irreducible curves on $X$ 
and $x\in X$ ranges over all the closed points of $X$.
The two boundary maps are given respectively by tame symbols for $K_2$ and 
by divisors of rational functions. 
Similarly $CH^3(X,2)$ is the cohomology of the complex
$$ K_3(k(X)) \to \underset{Z\subset X}\to{\bigoplus} K_2(k(Z)) \to 
 \underset{x\in X}\to{\bigoplus} k(x)^*,$$
where the second boundary map is given by tame symbols and the first by
localization theory for algebraic $K$-theory.
We note that, by a result of Merkurjev and Suslin, 
one is allowed to replace the Quillen $K_3$ by the Milnor $K^M_3$
and then the first boundary map is given also by tame symbols for $K^M_3$.
\vskip 5pt

For the study of the above groups we fix
$X\supset Z=\cup_{i\in I} Z_i\sp (I=\{1,2,\dots,m\})$, 
a simple normal crossing divisor on $X$ and consider particularly higher 
cycles supported on $Z$: we write
$$ CH^1(Z,1)=\Ker(\underset{i\in I}\to{\bigoplus} k(Z_i)^* \to 
 \underset{x\in Z}\to{\bigoplus} \Bbb Z),$$
$$ CH^2(Z,2)=\Ker(\underset{i\in I}\to{\bigoplus} K_2(k(Z_i)) \to 
 \underset{x\in Z}\to{\bigoplus} k(x)^*).$$
These are the Bloch's higher Chow groups of $Z$ (cf. [Bl]) and we have the 
exact sequence for $r=1,2$
$$ CH^{r+1}(U,r+1) \to CH^r(Z,r) \to CH^{r+1}(X,r).$$
Our first question is if one can find any interesting elements in $CH^r(Z,r)$
whose images in $CH^{r+1}(X,r)$ or $CH^{r+1}(X,r)^{ind}$,
the so-called indecomposable part of it,
are non-torsion. Our first main result Th.(0-1) suggests that this is impossible 
if $X\subset \bP^3$ is a very general hypersurface of sufficiently high degree
and the components of $Z$ are very general hypersurface sections of $X$ 
(see Def.(1-1) for the definition of $(X,Z)$ being very general). 
In order to state them, we need to introduce the indecomposable parts of
$CH^r(Z,r)$ :
$$ CH^r(Z,r)^{ind}=
\Coker\big(\underset{i\in I}\to{\bigoplus} CH^r(Z_i,r) \to CH^r(Z,r)\big).
$$
For an alternative description of $ CH^r(Z,r)^{ind}$ see Pr.(1-1).

\Th 0-1. \it Let $X\subset \bP^3$ be a very general hypersurface of degree 
$d$ and let $Z=\cup_{i\in I} Z_i$ with $Z_i\subset X$, a very general 
hypersurface sections of degree $e_i$. 
\roster
\item
If $d\geq 5$, $CH^2(U,2) \to CH^1(Z,1)^{ind}$ is surjective. 
\item
Assume $d\geq 6$ and that $(e_i,e_j,e_l)\not=(1,1,2)$ for distinct 
$i,j,l\in I$. Then
$CH^3(U,3) \to CH^2(Z,2)^{ind}$ is surjective. 
\endroster
\vskip 5pt\rm

The proof is easily reduced to the case where our base field $k=\Bbb C$.
Then the essential idea of the proof goes back to Griffiths' fundamental work
on algebraic cycles (cf. \cite{Gri}) and new improvements made later by
Green \cite{G1} and Voisin \cite{V}. The assumption on the generality allows us to extend 
our varieties to a family $(X,Z)/S$ of varieties over a large parameter space 
$S$ with the fibers $(\Xs,\Zs)$ over $s\in S$.
Then, by using the theory of variation of Hodge structures, we construct the
cycle class map for $CH^r(Z,r)$ 
$$ \pZSr: CH^r(Z,r) \to H^0(S,\PZr),$$
where $\PZr$ is a sheaf on $S_{an}$, the analytic site on $S(\Bbb C)$.
We have
$$\PZr =
\left\{\eqalign{
&\HQZX {3} 2 \cap F^{2}\HOZX {3},\quad \hbox{ if } r=1, \cr
\Ker\big(&\HOZX {3}/\HQZX {3} 3 @>\nabla>>\WS {1}\otimes 
\HOZX {3}/F^{2}\HOZX {3}\big),\quad \hbox{ if } r=2, \cr
}\right.$$
where $\HQZX q n$ is the local system on $S_{an}$ with fibers 
$H^q_{\Zs}(\Xs,\Bbb Q(n))$, the local cohomology of $\Xs$ with support in 
$\Zs$, and $\nabla:\HOZX q\to \WS 1\otimes \HOZX q$ is the sheaf of 
holomorphic sections of $\HCZX q=\HQZX q n\ottC$ with the Gauss-Manin 
connection and $F^p\HOZX q\subset \HOZX q$ is the holomorphic subbundle whose 
fibers give the Hodge filtration on $H^q_{\Zs}(\Xs,\Bbb C)$ defined by 
Deligne \cite{D2}. The above map is essentially given rise to by the regulator map 
from higher Chow groups to Deligne cohomology (cf. \cite{EV}).
For $r=2$ it is an analogue of the Abel-Jacobi map defined by Griffiths
\cite{Gri} after which we call the sections of $\PZr$ normal functions.
Now the key to the proof of Th.(0-1) is that one can compute the space
$H^0(S,\PZr)$ by using the theory of generalized Jacobian rings developed
in \cite{AS1}. A preliminary version of this method has been used in
\cite{SMS} to prove the vanishing of Deligne classes for $d \ge 5$ in the
case $r=1$. The key steps of the present computation have already been
carried out in \cite{AS2} and \cite{AS3}, where it has been applied to the so-called 
\it Beilinson's Hodge and Tate conjecture \rm for $U=X-Z$.
The disturbing assumption on the $e_i$'s in Th.(0-1)(2) is caused by a technical
obstruction in the Jacobian ring computation in \cite{AS3}.
\vskip 6pt

The second objective of the paper is to apply our Hodge theoretic 
invariants for the purpose of detecting non-trivial elements in 
the indecomposable part $CH^{r+1}(X,r)^{ind}$ of $CH^{r+1}(X,r)$. 
As Th.(0-1) suggests it is not hopeful for very general $(X,Z)$, while one may 
still hope for the possibility to find non-trivial examples 
among either special families or complete families of surfaces of low degree. 
This is done in \S5,\S6 and \S7.
After presenting the necessary formalism of infinitesimal invariants
of normal functions in \S4, we will show the following results:

\Th 0-2. \it Consider the family 
$$ X_{u,v}= \{F_{u,v}=x_0^5 + x_1x_2^4+x_2x_1^4+x_3^5 + ux_1^2x_2^3
+vx_0x_3K(x_0,...,x_3)=0 \}, \quad u,v \in {\Bbb C}$$ 
of quintic surfaces over $\Spec(\Bbb C[u,v])$ where
$ K$ is a homogenous polynomial of degree 3 with coefficients in 
$\Bbb C[u,v]$. Then there exist elements $\alpha_{u,v}$ in $CH^3(X_{u,v},2)$ 
supported on $Z=X\cap \{x_0x_3=0\}$
such that, for $u,v \in {\Bbb C}$ and $K$ very general, 
these elements are indecomposable modulo the image of $Pic(X_{u,v}) \otimes K_2(\bC)$.
\rm \vskip 6pt

\Th 0-3. \it On the family 
$$ X_{u}= \{ (x_0:...:x_3) \in {\Bbb P}^3 
\mid F_u(x)=x_0x_1^4+ x_1x_2^4+x_2x_0^4+ x_3^5+ ux_3x_1^4=0 \}, 
\quad u \in {\Bbb C}
$$ 
of quintic surfaces, there exist elements $\alpha_u$ in $CH^2(X_{u},1)$ 
supported on $Z=X\cap \{x_0x_3=0\}$ such that, for $u$ very general, 
these elements are indecomposable modulo the image of $Pic(X_{u}) \otimes 
\bC^*$.
\rm \vskip 6pt

The following examples were provided to us by Alberto Collino in a letter
from September 19, 1999. We are very grateful to him for letting us reproduce
the contents here. His result shows in particular that in case of surfaces of 
low degree, even a very general surface can carry indecomposable
cycles in $CH^3(S,2)$ which is supported on a smooth hyperplane section and 
which need not be rigid on the surface $S$:

\Th 0-4. \it On every very general quartic $K3$-surface $S$, there exists a
1-dimensional family of elements $\alpha_t$ in $CH^3(S,2)$ 
supported on a smooth hyperplane section of $X$ such that, for $t$ very 
general,
these elements are indecomposable modulo the image of $Pic(S)\otimes K_2(\bC)$.
\rm

\vskip 20pt

\head \S1. Higher Chow groups of normal crossing divisors. \endhead
\vskip 8pt 

In this section we recall some basic facts on Bloch's higher Chow groups 
and state main results from which Th.(0-1) is deduced.
In the whole paper we consistently neglect torsion and let $M$ 
denote $M\otZQ$ for an abelian group $M$. In what follows we fix the following.
\roster
\item"$(i)$"
$S=\Spec R$ is an affine smooth scheme over a field of characteristic zero.
\item"$(ii)$"
$X\hookrightarrow \bP^3_S=\Proj R[X_0,X_1,X_2,X_3]$ is a family of smooth 
hypersurfaces of degree $d$ with $f:X\to S$ the natural map.
\item"$(iii)$"
$Z_i\hookrightarrow X$ for $i\in I=\{1,2,\dots,m\}$ is a family of 
smooth hypersurface sections defined by a homogeneous polynomial 
$G_i\in R[X_0,X_1,X_2,X_3]$ of degree $e_i$ such that 
$Z=\cup_{i\in I} Z_i\subset X$ is a relatively normal crossing divisor.
We write for $1\leq j\leq m-1$
$$g_j=(G_j)^{e_s}/(G_s)^{e_j} \in \Cal O_{Zar}(U)^*=CH^1(U,1).$$
\endroster
\vskip 6pt

Fix integers $r\geq 1$.
The objects of our study are the higher Chow groups $CH^{r}(Z,r)$ (cf. \cite{Bl}),
particularly in the case $r=1,2$. Define
$$ \ZZ 1=\underset{1\leq i \leq m}\to{\coprod} Z_{i}\qaq
\ZZ 2 =\underset{1\leq i<j \leq m}\to{\coprod} Z_{i}\cap Z_{j}.$$

\Pr 1-1. \it We have the exact sequence
$$ CH^{r}(\ZZ 1,r)\to CH^{r}(Z,r) @>{\pi_Z}>> CH^{r}(Z,r)^{ind}\to 0,$$
where 
$CH^{r}(Z,r)^{ind}=\Ker\big(CH^{r-1}(\ZZ 2,r-1)\to CH^{r}(\ZZ 1,r-1)\big).$
\rm\pr
This follows immediately from the spectral sequence
$$ \sE a b 1= CH^{a+n}(\ZZ {-a},1-b)\Rightarrow CH^{n-1}(Z,-a-b).$$
\qed
\vskip 6pt

Over $\Spec(\bQ)$, we consider
$$ \tilde{M} \subset 
\overset{\vee}\to{\Bbb P}(H^0(\Bbb P^3,\Cal O_{\Bbb P^3}(d)))\times 
\prod_{i\in I} 
\overset{\vee}\to{\Bbb P}(H^0(\Bbb P^3,\Cal O_{\Bbb P^3}(e_i))),$$
the moduli space of $(X\supset Z=\cup_{i\in I} Z_i)$ where
$X\hookrightarrow \bP^3$ is a smooth hypersurface 
of degree $d$ and $Z_i\hookrightarrow X$ for $i\in I$ is a family of 
smooth hypersurface sections of degree $e_i$ such that 
$Z=\cup_{i\in I} Z_i\subset X$ is a normal crossing divisor.
The algebraic group $G=\PGL_{4}$ acts naturally on $\tilde{M}$.
By \cite{GIT}, there exists a dense open subset 
$\tilde{M}'\subset \tilde{M}$ 
such that
\roster
\item
$\tM'$ is stable under the action of $G$,
\item
the geometric quotient $M=\tilde{M}'/G$ exists and it is smooth over
$\Spec(\bQ)$.
\item
the universal family over $\tilde{M}$ descends to the family
$(\calX\hookleftarrow \calZ =\cup_{i\in I} \calZ_i)/M$.
\endroster
\vskip 6pt

\Def 1-1. \it Let $k$ be a field of characteristic zero.
\roster
\item
A family $(X,Z)/S$ over a field $k$, such as in the beginning of this section, is 
complete, if there exist a dominant, rational map $\pi:S\to M$ of schemes such that 
$(X\hookleftarrow Z)/S$ is the pullback of $(\calX\hookleftarrow \calZ)$ via $\pi$. 
By the same way, in the case $I=\emptyset$, we define a family $X/S$ 
of hypersurfaces of degree $d$ to be complete.
\item
A pair $(X,Z)$ of a smooth hypersurface $X\subset \bP^3$ and a normal crossing 
divisor $Z\subset X$ defined over $k$ is very general, if it is complete as 
a family over $S=\Spec(k)$.
\endroster
\rm\vskip 6pt

Now Th.(0-1) is a direct consequence of the following:

\Th 1-1. \it If $(X,Z)/S$ is complete and $d\geq 5$, 
$CH^2(U,2) \to CH^1(Z,1)^{ind}$ is surjective.
\vskip 5pt\rm

\Th 1-2. \it Assume $(X,Z)/S$ is complete and $d\geq 6$ and that 
$(e_i,e_j,e_l)\not=(1,1,2)$ for distinct $i,j,l\in I$. Then
$CH^3(U,3) \to CH^2(Z,2)^{ind}$ is surjective. 
\vskip 5pt\rm

\vskip 20pt

\head \S2. Normal functions associated to higher cycles on $Z$. \endhead

\vskip 8pt 

This section contains preliminary technical results for the proofs of Th.(1-1) and Th.(1-2)
in the next section. First we note that, by a well-known argument, we can reduce
our problem to the case where our base field $k$ is equal to $\Bbb C$, the
fields of complex numbers. Let $S,X,Z$ be as in the beginning of \S1 and assume that 
$S$ is a smooth affine variety over $\bC$. 
Let $j:U=X-Z \hookrightarrow X$ and
$i:Z \hookrightarrow X$ be the inclusions.
Consider the following complexes
of sheaves on $\Xan$ introduced by Deligne (cf. \cite{EV}): 
$$ \eqalign{
&\QDX n =\Cone(\bQ(n)\oplus \WX {\geq n} \to \WX \cdot)[-1],\cr
&\QDU r =\Cone(\bR j_* \bQ(r) \oplus \WXlogZ {\geq n}
\to \Bbb R j_*\Omega_U^\cdot)[-1],\cr
&\QDZX n =\Cone(\QDX n \to \QDU n)[-1].\cr
}$$
By \cite{D1}[Pr.3.1.8] we have the following distinguished triangles in 
$D^b(X_{an})$, the derived category 
of complexes of sheaves of bounded constructible cohomologies,
$$ i_* \bR i^! \bQ(n)_X \to \WXlogZ {<n}/W_0 [-1] \to \QDZX n \map d,
\leqno(2-1)$$
where $\WXlogZ {<n}$ is the truncated de Rham complex of $X$ with 
logarithmic poles along $Z$ and 
\par\noindent
$W_p=W_p\WXlogZ {\cdot}\subset \WXlogZ {\cdot}$ denotes the weight 
filtration. We define the following sheaves on $S_{an}$ 
$$ \cHDXn q=R^q f_* \QDX n ,\quad \cHDUn q=R^q f_* \QDU n,
\quad \cHDZXn q=R^q f_* \QDZX n.$$
Our fundamental tool to study higher Chow groups is the following commutative 
diagram
$$\matrix
CH^{n}(U,r+1) &@>{\oU{n}{r+1}}>>& H^0(S,\cHDUn {2n-r-1}) \\
\downarrow&&\downarrow\\
CH^{n-1}(Z,r) &@>{\oZ{n}{r}}>>& H^0(S,\cHDZXn {2n-r}) \\
\downarrow&&\downarrow\\
CH^{n}(X,r) &@>{\oX{n}{r}}>>& H^0(S,\cHDXn {2n-r}) \\
\endmatrix \leqno(2-2)$$
with the right vertical sequence arising from the following 
distinguished triangle in $D^b(X_{an})$
$$ \QDZX n \to \QDX n \to \QDU n \to $$
and the left vertical arrow arising from localization theory for higher Chow groups. 
The horizontal maps are defined by the theory of Chern classes in Deligne cohomology. 
The commutativity of the diagram is a consequence of the functoriality 
of Chern class maps into Deligne cohomology (cf. \cite{J}).
\vskip 6pt

In order to go further, we must introduce some more notations. 
For integers $q,n$ put
$$ \eqalign{
&\HQX q n=R^q f_* \bQ(n),\cr
&\HQU q n=R^q f_* \bR j_*\bQ(n),\cr
&\HQZX q n=R^q f_* i_*\bR i^!_*\bQ(n).\cr}$$
These are local systems on $\San$ and we have the long exact sequence
$$ \cdots\to \HQU {q-1} n \to \HQZX q n\to \HQX q n \to \HQU q n \to \cdots.$$
We also consider the local systems
$$ \HCX q=\HQX q n \otimes_{\bQ}\bC,
\quad \HCU q=\HQU q n \otimes_{\bQ}\bC,
\quad \HCZX q =\HQZX q n \otimes_{\bQ}\bC$$
and locally free $\cO$-modules 
$$ \HOX q=\HCX q\otimes_{\bC}\cO,
\quad \HOU q=\HCU q\otimes_{\bC}\cO,
\quad \HOZX q=\HCZX q\otimes_{\bC}\cO,$$
where $\cO$ denotes the sheaf of holomorphic functions on $S$.
By the relative Poincar\'e lemma we have the canonical isomorphisms
$$ \eqalign{
&\HOX q \isom R^q f_* \WXS{\cdot},\cr
&\HOU q \isom R^q f_* \WXSZ{\cdot},\cr
&\HOZX q \isom R^{q-1} f_* (\WXSZ{\cdot}/W_0),\cr}$$
where $\WXSZ \cdot$ is the relative de Rham complex with logarithmic poles
along $Z$ and $W_p\subset \WXSZ \cdot$ is the weight filtration (cf. \cite{D1}).
It follows from Deligne's theory of mixed Hodge structures that they are locally 
free $\cO$-modules and that the natural maps
$$ \eqalign{
 &F^n\HOX {q}:=R^qf_* \WXS {\geq n} \to  \HOX q,\cr  
 &F^n\HOU {q}:=R^qf_* \WXSZ {\geq n} \to  \HOU q,\cr  
 &F^n\HOZX {q}:=R^{q-1} f_* (\WXSZ {\geq n}/W_0) \to  \HOZX q\cr}$$  
are injective and their images are locally free $\cO$-submodules.
We have a long exact sequence
$$ \cdots\to F^n\HOZX q\to F^n\HOX q\to F^n\HOU q\to F^n\HOZX {q+1}\to\cdots
\leqno(2-3)$$
and we have integrable connections
$$ \eqalign{
&\nabla_X\scs \HOX {q} \to \WS 1\otimes \HOX {q},\cr
&\nabla_U\scs \HOU {q} \to \WS 1\otimes \HOU {q},\cr
&\nabla_Z\scs \HOZX {q} \to \WS 1\otimes \HOZX {q},\cr}
\leqno(2-4)$$
such that 
$ \Ker(\nabla_X)=\HCX q$, $\Ker(\nabla_U)=\HCU q$ and $\Ker(\nabla_Z)=\HCZX q$.
They satisfy
$$ \eqalign{
&\nabla_X(F^n \HOX {q})\subset \WS 1\otimes F^{n-1}\HOX {q},\cr
&\nabla_U(F^n \HOU {q})\subset \WS 1\otimes F^{n-1}\HOU {q},\cr
&\nabla_Z(F^n \HOZX {q})\subset \WS 1\otimes F^{n-1}\HOZX {q}.\cr
}$$
We now define the \it sheaf of normal functions \rm with support in $Z$:

\Def 2-1. \it Assume $n=r+1$ and $r=1$ or $r=2$. 
We define
$$\PZr =
\left\{\eqalign{
\Ker\big(&\HOZX {3}/\HQZX {3} 3 @>\nabla>>\WS {1}\otimes 
\HOZX {3}/F^{2}\HOZX {3}\big),\quad \hbox{ if } r=2, \cr
&\HQZX {3} 2 \cap F^{2}\HOZX {3},\quad \hbox{ if } r=1, \cr
}\right.$$
$$\PUr =
\left\{\eqalign{
\Ker\big(&\HOU {2}/\HQU {2} 3 @>\nabla_U>>\WS {1}\otimes 
\HOU {2}/F^{2}\HOU {2}\big),\quad \hbox{ if } r=2, \cr
&\HQU {2} 2 \cap F^{2}\HOU {2},\quad \hbox{ if } r=1, \cr
}\right.$$
\rm

\Pr 2-1. \it For $r=1$ and $r=2$, diagram (2-2) induces the commutative diagram
$$\matrix
CH^{r+1}(U,r+1) &@>{\pUSr}>>& H^0(S,\PUr)\\
\downarrow&&\downarrow\\
CH^{r}(Z,r) &@>{\pZSr}>>& H^0(S,\PZr),\\
\endmatrix$$
where the right vertical arrow is induced by the localization map
$\HQU{q}{n} \to \HQZX {q+1}{n}$.

\rm\demo{Proof}
We filter $\WXlogZ {< n}/W_0$ by the subcomplexes
$$ F^{p}_S (\WXlogZ {< n}/W_0)=
\Im(f^*\WS p\otimes \WXlogZ {< n-p}\to \WXlogZ {< n}/W_0)$$
so that we have
$$ Gr^p_{F_S}(\WXlogZ {< n}/W_0)=f^*\WS p\otimes (\WXSZ {<n-p}/W_0)[-p].$$
The projection formula now gives
$$ R^{p+q}f_* Gr^p_{F_S}(\WXlogZ {< n}/W_0) = 
\WS p\otimes \HOZX q/F^{n-p}.$$
In view of the exact sequence (2-1), this gives rise to 
the spectral sequence 
$$ \sIE pq1(n)\Rightarrow \cHDZXn {p+q} \quad \hbox{ with }
\sIE pq1(n)=
\left\{\eqalign{
& \HQZX q n \quad \hbox{ if } p=0, \cr
& \WS {p-1} \otimes \HOZX q/F^{n-p+1}\quad \hbox{ if } p\geq 1. \cr
}\right.
\leqno(2-5)$$
Note that the $d_1$-differential of the spectral sequence is induced by the
connection $\nab Z$. In view of Lem.(2-1) below we have the edge homomorphism
$\cHDZXn {r+2} \to \PZr$ (Note $\PZr=\sIE {r-1} {3} 2(r+1)$).
Now $\pZSr$ is defined to be the composite of this map and $\oZ{r+1}{r}$ 
(cf. (2-2)). 
The construction of $\pUSr$ is parallel to that of $\pZSr$ by using 
the following spectral sequence
$$ \sE pq1(n)\Rightarrow \cHDUn {p+q} \quad \hbox{ with }
\sE pq1(n)=
\left\{\eqalign{
& \HQU q n \quad \hbox{ if } p=0, \cr
& \WS {p-1} \otimes \HOU q/F^{n-p+1}\quad \hbox{ if } p\geq 1. \cr
}\right.$$
The compatibility with localization sequences and with Chern class
maps to Betti cohomology follows from basic properties of Chern class
maps in Deligne cohomology (cf. \cite{EV, J}).
\qed
\enddemo
\vskip 6pt

Here are some basic vanishing results we need later:

\Lem 2-1. \it 
$\sIE p q 1(n) =0$ either if $q\leq 1$ or $p\not=1,p+q\leq n$ or 
$p\not=1,q=n-p+1>3$.
\rm\pr
The lemma follows from Lem.(2-2) below in view of the exactness of the sequence
$$ 0\to \HCZX q\to \HOZX q @>{\nab Z}>> \WS 1\otimes\HOZX q @>{\nab Z}>>
\WS 1\otimes\HOZX q @>{\nab Z}>> \cdots.$$
\vskip 5pt

\Lem 2-2. 
\roster
\item
$\HOZX 1=0$.
\item
$F^m\HOZX {q+1}=0$ for $m\geq q+1$.
\item
$F^q\HOZX {q+1}=0$ for $q>2=\dim(X/S)$.
\endroster
\pr
For $p\geq 0$ we denote
$$ \ZZ p=
\left\{\eqalign{
&X \quad \hbox{ if } p=0, \cr
&\underset{1\leq i \leq m}\to{\coprod}Z_{i} \quad \hbox{ if } p=1, \cr
&\underset{1\leq i<j \leq m}\to{\coprod}Z_{i}\cap Z_{j}
 \quad \hbox{ if } p=2, \cr
&\emptyset \quad \hbox{ if } p\geq 3. \cr
}\right.$$
We have the spectral sequence
$$\sE ab1 =F^{r+a}H^{b+2a}_{\cO}(\ZZ {-a}/S)\Rightarrow F^r\HOZX {a+b+1},
\leqno(2-6)$$
which is constructed from the weight filtration $W_p\subset \WXSZgeq r$
with the isomorphism
$$ Gr^W_p \WXSZgeq r \isom (i_p)_*\Omega_{\ZZ p/S}^{\geq r-p}[-p],$$
where $i_p:\ZZ p \to X$ is the natural morphism. 
The lemma follows from this noting that
$$ F^{r+a}H^{q+a}_{\cO}(\ZZ {-a}/S)=0 \qif r>\max\{q,d\},$$
since $\dim(\ZZ {-a}/S)=\dim(X/S)+a$.
\qed
\vskip 6pt

\Pr 2-2. \it Let $x\in CH^r(Z,r)$ and assume $\pZSr(x)=0$. Then $\pi_Z(x)=0$,
if $r \le 2$ and where $\pi_Z$ is as in Pr.(1-1).
\rm\pr 
By localization theory we have the maps
$$ \eqalign{
&\HQZX 3 2 \to H^0_{\bQ}(\ZZ 2/S)=g_*\bQ,\cr
&\HOZX 3/\HQZX 3 3 \to H^0_{\cO}(\ZZ 2/S)/H^0_{\bQ}(\ZZ 2/S)(1),
=g_*(\cO_{\ZZ 2}/\bQ(1)),\cr
}$$
where $g:\ZZ 2\to S$ is the natural morphisms. 
It induces the commutative diagrams
$$\matrix
CH^1(Z,1)&@>{\pZS 1}>>&H^0(S,\PZ 1)\\
\downarrow\rlap{$\pi_Z$}&&\downarrow\\
CH^{0}(\ZZ 2)&@>{\alpha}>>&H^0(\ZZ 2,\bQ)\\
\endmatrix
\quad\hbox{ and }\quad
\matrix
CH^2(Z,2)&@>{\pZS 2}>>&H^0(S,\PZ 2)\\
\downarrow\rlap{$\pi_Z$}&&\downarrow\\
CH^{1}(\ZZ 2,1)&@>{\beta}>>&H^0(\ZZ 2,\cO_{\ZZ 2}/\bQ(1)),\\
\endmatrix$$
where $\alpha$ is the obvious map and $\beta$ is the isomorphism
$CH^{1}(\ZZ 2,1)\isom \Cal O_{Zar}(\ZZ 2)^*$ followed by the logarithm. 
Either of the maps is injective and the desired assertion is proven in view
of Pr.(1-1). 
\qed

\vskip 20pt

\head \S3. Proof of the main results. \endhead
\vskip 8pt 

In this section we prove Th.(1-1) and Th.(1-2).
By Pr.(2-2) they follow from the following claims.
Recall that $\delta: CH^{r+1}(U,r+1) \to CH^r(Z,r)$ denotes the boundary map.

\Claim 3-1. \it Under the assumption of Th.(1-1),
$\Im(\pZS 1)=\Im(\pZS 1 \cdot\delta)$.
\rm\vskip 5pt

\Claim 3-2. \it Under the assumption of Th.(1-2),
$\Im(\pZS 2)=\Im(\pZS 2 \cdot\delta)$.
\rm\vskip 5pt

First we prove Claim(3-1).

\Def 3-1. \it With the notation as in the beginning of \S1, 
let $\MQU\subset \HQU 2 2$ be the (constant) subsheaf generated
by the sections $\omega_{ij}:=\pUS 1(\{g_i,g_j\})$ for $1\leq i<j\leq m-1$.
\rm\vskip 4pt

We need the following result from [AS2].

\Th 3-1. \it If $(X,Z)$ is complete and $d\geq 4$, $H^0(S,\HQZX 3 2)=\MQU$.
In particular $\pUS 1$ is surjective.
\rm\vskip 5pt

By the above theorem, Claim(3-1) in the case $r=1$ follows from the following.

\Pr 3-1. \it Assuming $(X,Z)$ is complete and $d\geq 5$,
$H^0(S,\HQU 2 2)\to  H^0(S,\HQZX3 2)$ is surjective.
\rm\vskip 5pt

We note that it suffices to show the surjectivity of 
$H^0(S,\HCU 2)\to  H^0(S,\HCZX3)$.

\Lem 3-1. \it Under the assumption of Pr.(3-1),
$H^0(S,\HQZX3 2)\subset H^0(S,\HCZX3\cap F^2\HOZX 3).$
\rm\demo{Proof} 
Localization theory provides the exact sequence
$0\to \KQ \to \HQZX 3 2 \to \JQ \to 0$
with 
$$ \KQ=\bigoplus_{i\in I} \HQZi 1 1,\quad
\JQ=\Ker(\HQ^0(\ZZ 2/S)\to\bigoplus_{i\in I} \HQZi 2 1)).$$
We have $H^0(S,\KQ)=0$, which follows from [AS2, Th.(4-2)]
or the monodromy argument of Green and Voisin (cf. [G1]).
This proves the desired assertion.
\qed
\enddemo
\vskip 5pt

By Lem.(3-1) the proof of Pr.(3-1) is reduced to the following.

\Lem 3-2. \it Under the assumption of Pr.(3-1), the following map is surjective
$$H^0(S,\HCU 2\cap F^2\HOU 2)\to H^0(S,\HCZX3\cap F^2\HOZX 3).$$
\rm\demo{Proof} 
Note that we have the exact sequence
$$ 0\to \HQXpr 2 2 \to \HQU 2 2 \to \HQZX 3 2 \to 0,$$
where $\HQXpr 2 2$ is the primitive part of $\HQX 2 2$.
It induces the following commutative diagram
$$\matrix
0&&0&&0\\
\downarrow&&\downarrow&&\downarrow\\
F^2\HOXpr 2 &@>{\nabla}>>& \WS 1\otimes F^1\HOXpr 2 &@>{\nabla}>>& 
\WS 2\otimes F^0\HOXpr 2 \\
\downarrow&&\downarrow&&\downarrow\\
F^2\HOU 2 &@>{\nabla}>>& \WS 1\otimes F^1\HOU 2 &@>{\nabla}>>& 
\WS 2\otimes F^0\HOU 2 \\
\downarrow&&\downarrow&&\downarrow\\
F^2\HOZX 3 &@>{\nabla}>>& \WS 1\otimes F^1\HOZX 3 &@>{\nabla}>>& 
\WS 2\otimes F^0\HOZX 3 \\
\downarrow&&\downarrow&&\downarrow\\
0&&0&&0,\\
\endmatrix$$
where the vertical sequences are exact. Noting that $S$ is affine, it implies 
that Lem.(3-2) follows from the exactness at the middle of the upper 
horizontal sequence, which in turn follows from the following result, a 
consequence of the theory of Jacobian rings (cf. [G2]).
\enddemo

\Pr 3-2. \it Consider the following complex 
$$ \WS {p-1}\otimes \HHOXpr {1+q}{1-q} @>{\nabb X}>>
 \WS {p}\otimes \HHOXpr {q} {2-q} @>{\nabb X}>>
 \WS {p+1}\otimes \HHOXpr {q-1}{3-q},$$
where 
$\HHOX a b=Gr^a_F\HOX {a+b}=R^b f_*\WXS a$
and $\HHOXpr ab$ denotes the primitive part of $\HHOX ab$.
By convention $\WS r=0$ if $r<0$.
Assuming $(X,Z)$ is complete, the complex is exact if 
$q\geq 1$ and $dq\geq p+4$.
\rm\vskip 6pt

We now show Claim(3-2).

\Def 3-2. \it We define
$$\FZ 2 =
\Ker\big(\WS {1}\otimes F^{2}\HOZX {3} @>\nabla>>
\WS {2}\otimes F^{1}\HOZX {3}\big),$$
$$\FU 2 =
\Ker\big(\WS {1}\otimes F^{2}\HOU {2} @>\nabla>>
\WS {2}\otimes F^{1}\HOU {2}\big).$$
By definition we have the exact sequences
$$ 0\to \HCZX 3/\HQZX 3 3\to \PZ 2 @>{\nabla}>> \FZ 2,$$
$$ 0\to \HCU 2/\HQU 23\to \PU 2 @>{\nabla}>> \FU 2.$$
We have the following maps induced by $\pZS 2$ and $\pUS 2$ respectively
$$\cZS 2 :CH^2(Z,2)\ottC\to H^0(S,\FZ 2),\quad
\cUS 2 :CH^3(U,3)\ottC\to H^0(S,\FU 2).$$
\rm\vskip 5pt

\Def 3-3. \it With the notation as in the beginning of \S1, let 
$$\NCU\subset H^0(S,\FU 2)\subset H^0(S,\WS {1}\otimes F^{2}\HOU {2})$$
be the $\Bbb C$-linear subspace generated by 
$\cUS 2(\{g_i,g_j,g_k\})$ with $1\leq i<j<k\leq m-1$.
\rm\vskip 4pt

We need the following result from [AS3].

\Th 3-2. \it Under the assumption of Th.(1-2), 
$$\WSd 1\otimes\MCU \oplus \NCU \isom \FU 2,$$ 
where $\MCU=\MQU\ottC\subset \HCU 2\cap F^2\HOU 2.$ 
\rm\vskip 5pt

Thanks to Th.(3-2), Claim(3-2) follows from the following three lemmas.

\Lem 3-3. \it Consider the composite map
$$ \alpha:\NCU \oplus \MCU\otimes H^0(S,\WSd 1) \to H^0(S,\FU 2)\to 
H^0(S,\FZ 2).$$
Assume 
$\alpha(x+y)\in \Im(\cZS 2)$ for $y\in \NCU$ and 
$$x=\sum_{1\leq i<j\leq m-1}\omega_{ij}\otimes c_{ij}\in 
\MCU\otimes H^0(S,\WSd 1)\quad \hbox{(cf. Def.(3-1))}.$$
Then $c_{ij}$ is contained in the image of 
$d\log:CH^1(S,1)\ottC \to H^0(S,\WSd 1).$
\rm\vskip 5pt

\Lem 3-4. \it Assuming $(X,Z)$ is complete and $d\geq 5$, we have
$$ \Ker(H^0(S,\PZ 2)\to H^0(S,\FZ 2))=
\sum_{\overset{x\in \Bbb C^*}\to{1\leq i<j\leq m-1}}
\pZS 2(\delta\{g_i,g_j,x\}).$$
\rm\vskip 5pt

\Lem 3-5. \it Assuming $(X,Z)$ is complete and $d\geq 6$, 
$H^0(S,\FU 2)\to H^0(S,\FZ 2)$ is surjective.
\rm\vskip 6pt

\demo{Proof of Lem.(3-3)}
We have the commutative diagram
$$\matrix
CH^3(U,3)\ottC &@>{\delta}>>& CH^2(Z,2)\ottC &\to& CH^1(\ZZ 2,1)\ottC\\
\downarrow\rlap{$\cUS 2$}&&\downarrow\rlap{$\cZS 2$}
&&\downarrow\rlap{$d\log$}\\
H^0(S,\FU 2)&\to& H^0(S,\FZ 2) &@>{\beta}>>& H^0(\ZZ 2,\Omega^1_{\ZZ 2,d=0}),\\
\endmatrix
\leqno(3-3-1)$$
where $\beta$ is induced by 
$Res_Z:\HQZX 3 2\to \HQ^0(\ZZ 2/S)=g_*\bQ_{\ZZ 2}$ with $g:\ZZ 2\to S$ 
the projection. 
By a standard norm argument, we may prove Lem.(3-3) after a finite
\'etale base change of $S$, so that we may assume $\ZZ 2$ is a disjoint
union of copies of $S$.
This implies that the horizontal maps $\iota$ in the following commutative 
diagram are isomorphisms
$$ \matrix
H^0(\ZZ 2,\Bbb Q)\otimes CH^1(S,1)\ottC &@>{\iota}>>& 
CH^1(\ZZ 2,1)\ottC\\
\downarrow\rlap{$d\log$}&&\downarrow\rlap{$d\log$}\\
H^0(\ZZ 2,\Bbb Q)\otimes H^0(S,\WSd 1) &@>{\iota}>>& 
H^0(\ZZ 2,\Omega^1_{\ZZ 2,d=0}).\\
\endmatrix
\leqno(3-3-2)$$
By (3-3-1), the assumption of Lem.(3-3) implies
$$\sum_{1\leq i<j\leq m-1}\iota(Res_U(\omega_{ij})\otimes c_{ij})\in 
\Im(d\log:CH^1(\ZZ 2,1)\ottC \to H^0(\ZZ 2,\Omega^1_{\ZZ 2,d=0})),$$
where $Res_U:\MQU\to H^0(\ZZ 2,\Bbb Q)$ is the composite of $\MQU\to \HQZX 3 2$
and $Res_Z$. This implies the desired assertion due to (3-3-2) and 
the linear independence of $Res_U(\omega_{ij})$ with $1\leq i<j\leq m-1$,
which is easily checked.
\qed
\enddemo
\vskip 5pt

\demo{Proof of Lem.(3-4)}
By Th.(3-1) and Pr.(3-1) we have 
$$\MQU\otimes\Bbb C/\Bbb Q(1)\isom H^0(S,\HCZX 3/\HQZX 33)
\isom \Ker(H^0(S,\PZ 2)\to H^0(S,\FZ 2)),$$
where the injectivity of the above map follows that of $Res_U$ in the proof of
Lem.(3-3) and the second isomorphism follows from the first exact sequence in 
Def.(3-2). This proves Lem.(3-4) since we have 
$$ \pZS 2(\{g_i,g_j,x\})=\omega_{ij}\otimes \log x
\in \MQU\otimes \Bbb C/\Bbb Q(1)\quad \hbox{for } x\in \Bbb C^*.$$
\qed
\enddemo
\vskip 5pt

\demo{Proof of Lem.(3-5)}
As in the proof of Lem.(3-2) we have the commutative diagram
$$\matrix
0&&0&&0\\
\downarrow&&\downarrow&&\downarrow\\
\WS 1\otimes F^2\HOXpr 2 &@>{\nabla}>>& \WS 2\otimes F^1\HOXpr 2 
&@>{\nabla}>>& \WS 3\otimes F^0\HOXpr 2 \\
\downarrow&&\downarrow&&\downarrow\\
\WS 1\otimes F^2\HOU 2 &@>{\nabla}>>& \WS 2\otimes F^1\HOU 2 &@>{\nabla}>>& 
\WS 3\otimes F^0\HOU 2 \\
\downarrow&&\downarrow&&\downarrow\\
\WS 1\otimes F^2\HOZX 3 &@>{\nabla}>>& \WS 2\otimes F^1\HOZX 3 &@>{\nabla}>>& 
\WS 3\otimes F^0\HOZX 3 \\
\downarrow&&\downarrow&&\downarrow\\
0&&0&&0,\\
\endmatrix$$
where the vertical sequences are exact. From this Lem.(3-5) follows from 
the exactness at the middle of the upper horizontal sequence, which follows 
from Pr.(3-2).
\qed
\enddemo

\vskip 20pt

\head \S4. Indecomposable parts of infinitesimal invariants. \endhead

\vskip 8pt 
Let the notation be as in the beginning of \S1 and assume additionally the condition:
\roster
\item"$(4-1)$" $\HOX 1=0$.
\endroster

We want to capture elements of $CH^r(Z,r)$ whose image in 
$CH^{r+1}(X,r)$ are indecomposable in the sense defined below.
For this purpose we use the infinitesimal invariants of
normal functions. Recall the notation in \S2. 
By definition we have the natural maps
$$ \eqalign{
&\cHDXn q \to R^q f_* \WX {\geq n},\cr
&\cHDUn q \to R^q f_* \WXlogZ{\geq n} ,\cr
&\cHDZXn q \to R^{q-1} f_* \WXlogZ{\geq n}/W_0.\cr}.$$
We denote the above maps by the same letter $\delta$.
The commutative diagram (2-2) gives rise to the following commutative 
diagram
$$\matrix
CH^{n}(U,r+1) &@>{\delta\oU{n}{r+1}}>>& \Q(S,R^{2n-r-1}f_*\WXlogZ{\geq n})\\
\downarrow&&\downarrow\\
CH^{n-1}(Z,r) &@>{\delta\oZ{n}{r}}>>& \Q(S,R^{2n-r-1}f_*\WXlogZ{\geq n}/W_0) \\
\downarrow&&\downarrow\\
CH^{n}(X,r) &@>{\delta\oX{n}{r}}>>& \Q(S,R^{2n-r}f_*\WX{\geq n}) \\
\endmatrix \leqno(4-2)$$
with the right vertical sequence arising from the exact sequence
$$0\to \WX {\geq n}\to \WXlogZ {\geq n}\to \WXlogZ {\geq n}/W_0\to 0.$$
\vskip 6pt

We proceed in a similar manner for the construction of the
spectral sequence (2-5).
Filtering $\WXlogZ {\geq n}/W_0$ by the subcomplexes
$$ F^p_S \WXlogZ {\geq n}/W_0=
\Im(f^*\WS p\otimes \WXlogZ {\geq n-p}\to \WXlogZ {\geq n}/W_0),$$
we have
$$ Gr^p_{F_S}(\WXlogZ {\geq n}/W_0)=f^*\WS p\otimes (\WXSZgeq {n-p}/W_0).$$
Thus we obtain the spectral sequence
$$\sIIEZ p q 1= \WS p\otimes F^{n-p}\HOZX {q} \Rightarrow
R^{p+q-1}f_*\WXlogZ {\geq n}/W_0.\leqno(4-3)$$
Using the same arguments, we also construct the spectral sequence
$$\sIIEX p q 1= \WS p\otimes F^{n-p}\HOX {q} \Rightarrow
R^{p+q}f_*\WX {\geq n}.\leqno(4-4)$$
The following construction provides a more intrinsic definition of 
$\cZS 2$ given in Def.(3-2).

\Def 4-1. \it Let $n=r+1\geq 2$ be an integer.
\roster
\item
We define
$$\FZr =\Ker\big(\WS {r-1}\otimes F^2\HOZX {3}@>\nabla>>
\WS {r}\otimes F^{1}\HOZX {3}\big),$$
where $\nabla$ is induced by the connection (2-4) and coincides with 
$d_1$-differential of (4-3). By definition 
$\FZr=\sIIEZ {r-1}{3} 2$ and we define the map
$$ \cZSr\scs CH^r(Z,r) \to \Q(S,\FZr)$$
to be the composite of $\delta\oZ{r+1}{r}$ and the map induced by 
$R^{r+1}f_*\WXlogZ{\geq r+1}/W_0\to \FZr$ that is an edge homomorphism
of the spectral sequence (4-3).
\item
We define $\FXr=\sIIEX r 2 2$, which is the homology of the complex
$$\WS {r-1}\otimes F^2\HOX {2}@>\nabla>>
\WS {r}\otimes F^{1}\HOX {2}@>\nabla>>
\WS {r+1}\otimes F^{0}\HOX {2}$$
and define
$$ \cXSr \scs CH^{r+1}(X,r)\to \Q(S,\FXr)$$
to be the composite of $\delta\oX{r+1}{r}$ and the map induced by 
$R^{r+1}f_*\WX {\geq r+1} \to \FXr$ that is an edge homomorphism of (4-4).
\endroster
\vskip 6pt\rm

Note that the edge homomorphisms exist due to the fact that on one hand
$\sIIEZ p q 1=0$ if $q\leq 1$ or $p+q\leq n$ or $q=n-p+1>3$
by Lem.(2-2) and that on the other hand
$\sIIEX p q 1 =0$ unless $q=0,2,4$ by the assumption (4-1) and that
$\sIIEX p 4 1 =\WS p\otimes F^{r-p+1}\HOX 4=0$ if $p\leq r-2$.
\vskip 6pt

Next we consider the \it indecomposable part \rm
$$ CH^{r+1}(X,r)^{ind}=\Coker(CH^r(S,r)\otZ CH^1(X)\to CH^{r+1}(X,r)).$$
Write 
$$F^p\HNX q=\Ker(F^p\HOX q@>{\nabla_X}>> \WS 1\otimes F^{p-1}\HOX q),$$
which is $\bC$-vector space. Let
$$ cl_X\scs CH^1(X)\to \Q(S,F^1\HNX 2)$$
be the composite of the Chern class map
$\delta\oX 1 0:CH^1(X)\to \Q(S,R^2f_*\WX {\geq 1})$
and the map induced by
$$R^2f_*\WX{\geq 1}\to F^1\HNX 2$$
that is an edge homomorphism of (4-4) with $n=r+1=1$
($\sIIEX 0 2 2 =F^1\HNX 2$). 
Let
$$\psi_S^r\scs CH^r(S,r)\to \Q(S,\WSd r)$$
be induced by $\delta\oX r r$ with $X=S$ (cf. (4-2)). Define the map
$$ \alpha=\psi^r_S\otimes cl_X\scs 
CH^r(S,r)\otZ CH^1(X)\to \Q(S,\WSd r\otC F^1\HNX 2).$$ 

\Pr 4-1. \it
We have the commutative diagram
$$\matrix
CH^r(S,r)\otZ CH^1(X)&@>{\alpha}>>&\Q(S,\WSd r\otC F^1\HNX 2)\\
\downarrow\rlap{$\beta$}&&\downarrow\rlap{$\gamma$}\\
CH^{r+1}(X,r)&@>{\cXSr }>>& \Q(S,\FXr),\\
\endmatrix$$
where $\beta$ is induced by the product structure on higher Chow groups
and $\gamma$ is induced by the natural map
$$ \WSd r\otC F^1\HNX 2\to \WS r\otO F^1\HOX 2.$$
\rm\pr
We consider the following diagram
$$\matrix
CH^r(S,r)\otimes CH^1(X)&@>\beta>>& CH^{r+1}(X,r)\\
\downarrow\rlap{$\delta\oS r r\otimes\delta\oX 1 0$}&&
\downarrow\rlap{$\delta\oX {r+1} r$}\\
\Q(S,\WSd r)\otC \Q(S,R^2f_*\WX {\geq 1})&@>\epsilon>>& 
\Q(S,R^{r+1}f_*\WX {\geq r+1})\\
\downarrow&&\downarrow\\
\Q(S,\WSd r)\otC \Q(S,F^1\HNX 2)&@>{\gamma}>>& \Q(S,\FXr).\\
\endmatrix$$
Here $\oS r r$ and $\oX 1 0$ are Chern class maps (cf. (2-2))
and $\epsilon$ is induced by applying $R^{r+1}f_*$ to the product
$$ f^*\WS r\otimes \WX {\geq 1}[r-1]\to \WX {\geq r+1}.$$
The upper diagram is commutative due to the compatibility of Chern class map 
with product. The commutativity of the lower diagram is easily seen.
Now Pr.(4-1) follows from the fact that $\oS r r=\psi_S^r$. 
\qed
\vskip 6pt

\Def 4-2. \it
We define the indecomposable part of $\FXr$ 
$$ \FXindr=\Coker(\WSd r\otC F^1\HNX 2 \to \FXr),$$
and the map
$$ \cXSr \scs CH^{r+1}(X,r)^{ind}\to \Q(S,\FXindr),$$
which is induced by the commutative diagram of Pr.(4-1).
\rm
\vskip 6pt

Now the commutative diagram (2-2) gives rise to the commutative diagram
$$\matrix
 CH^{r}(Z,r)&@>{\cZSr}>>& \Q(S,\FZr)\\
\downarrow&&\downarrow\rlap{$\lambda$}\\
 CH^{r+1}(X,r)^{ind}&@>{\cXSr}>>& \Q(S,\FXindr),\\
\endmatrix\leqno(4-5)$$
where $\lambda$ is obtained as follows. We consider the commutative diagram
$$\matrix
\WS {r-1}\otimes F^2\HOX 2 &@>{\nabla_X}>>&
\WS {r}\otimes F^1\HOX 2 &@>{\nabla'_X}>>&\WS {r+1}\otimes F^0\HOX 2\\
\downarrow&&\downarrow\rlap{$\gamma$}&&\downarrow\\
\WS {r-1}\otimes F^2\HOU 2 &@>{\nabla_U}>>&
\WS {r}\otimes F^1\HOU 2 &@>{\nabla'_U}>>&\WS {r+1}\otimes F^0\HOU 2\\
\downarrow\rlap{$\beta$}&&\downarrow\\
\WS {r-1}\otimes F^2\HOZX 3 &@>{\nabla_Z}>>&\WS {r}\otimes F^1\HOZX 3\\
\downarrow\\0.\\
\endmatrix\leqno(4-6)$$
All vertical sequences are exact.
The surjectivity of $\beta$ is a consequence of the assumption (4-1). 
Putting
$$\HNZX q=\Ker(\HOZX q@>{\nabla_Z}>> \WS 1\otimes\HOZX q),$$
the spectral sequence (2-6) gives rise to the isomorphism
$$\HNZX 2\isom \bigoplus_{i\in I} \Ker(\cO_{Z_i}@>d>>\Omega_{Z_i}^1)
\isom \bigoplus_{i\in I} \bC.$$
Thus we have
$$\eqalign{\Ker(\gamma)
&=\WS r\otO \Im(\HOZX 2\to \HOX 2)\cr
&=\WS r\otC \Im(\HNZX 2\to \HNX 2)\cr
&\subset \WS r\otC F^1\HNX 2, \cr}$$
so that $\Ker(\gamma)\cap \Ker(\nabla'_X)\subset \WSd r\otC F^1\HNX 2$.
Hence the commutative diagram (4-6) gives rise to the map
$$ \Ker(\nabla_Z)\to 
\Coker\big( \WSd r\otC F^1\HNX 2\to \Ker(\nabla'_X)/\Im(\nabla_X)\big),$$
which is defined to be $\lambda$.
\vskip 6pt

For later use, we introduce the following linear version of the above 
construction.

\Def 4-3. \it
\roster
\item
We define $\CXr$ to be the homology of the complex of $\cO$-modules
$$\WS {r-1}\otimes \HHOX 2 0@>{\nabla_X}>>
\WS {r}\otimes \HHOXpr 1 1@>{\nabla_X}>>
\WS {r+1}\otimes \HHOX 0 2,$$
where $\HHOX p q =Gr^p_F \HOX {p+q}=R^qf_*\WXS p$ and 
$\HHOXpr 1 1$ is the primitive part of $\HHOX 1 1$.
We have the canonical map
$ \FXr \to \CXr$
and we put
$$ \CXindr=\Coker(\WSd r\otC F^1\HNX 2 \to \FXr\to \CXr).$$
The map $\cXSr$ induces 
$$ \tXSr \scs CH^{r+1}(X,r)^{ind}\to \Q(S,\CXindr).$$
\item
We define
$$\CZr =\Ker(\WS {r-1}\otimes \HHOZX 2 1@>{\nabla_Z}>>
\WS {r}\otimes \HHOZX 1 2),$$
where $\HHOZX p {q+1}=Gr^p_F\HOZX {p+q+1}=R^{p+q}f_*\WXSZ p/W_0$.
By Lem.(2-2)(1) we have the natural injection
$\FZr \hookrightarrow \CZr$ and the map $\cZSr$ induces 
$$ \tZSr \scs CH^r(Z,r)\to \Q(S,\CZr).$$
\rm
\endroster
\rm\vskip 5pt

The commutative diagram (4-5) gives rise to the commutative diagram
$$\matrix
 CH^{r}(Z,r)&@>{\tZSr}>>& \Q(S,\CZr)\\
\downarrow&&\downarrow\rlap{$\lambda$}\\
 CH^{r+1}(X,r)^{ind}&@>{\tXSr}>>& \Q(S,\CXindr),\\
\endmatrix\leqno(4-7)$$
where the right vertical arrow is induced by the commutative diagram
$$\matrix
\WS {r-1}\otimes \HHOX 2 0&@>{\nabla_X}>>&
\WS {r}\otimes \HHOXpr 1 1 &@>{\nabla_X}>>&\WS {r+1}\otimes \HHOX 0 2\\
\downarrow&&\downarrow&&\downarrow\\
\WS {r-1}\otimes \HHOU 2 0&@>{\nabla_U}>>&
\WS {r}\otimes \HHOU 1 1 &@>{\nabla_U}>>&\WS {r+1}\otimes \HHOU 0 2\\
\downarrow&&\downarrow\\
\WS {r-1}\otimes \HHOZX 2 1 &@>{\nabla_Z}>>&\WS {r}\otimes \HHOZX 1 2\\
\downarrow\\0\\
\endmatrix\leqno(4-8)$$
with $\HHOU p q =Gr^p_F \HOU {p+q}=R^qf_*\WXSZ p$.

\vskip 20pt

\def\Fuv{F_{u,v}} 
\def\Fu{F_{u}}
\def\sp{\hbox{ }}
\def\tRF{\tilde{R}_F}
\def\Spec{\hbox{\rm Spec}}
\def\Proj{\hbox{\rm Proj}}
\def\isom{@>\cong>>}
\def\Ker{\hbox{\rm Ker}}
\def\mod{\sp\hbox{\rm mod}\sp}

\head \S5. An Example in $CH^3(X,2)$. \endhead
\vskip 8pt 

We mostly keep the notations of the previous paragraphs.
In this section we prove Th.(0-2) which 
provides an example of classes in $CH^3(X,2)$ on a family of 
smooth projective complex surfaces $X_{u,v}$, which are indecomposable modulo
the image of $K_2({\Bbb C})\otimes CH^1(X)$. We state the result here again.

\Pr 5-1. \it Consider the family 
$$ X_{u,v}= \{F_{u,v}=x_0^5 + x_1x_2^4+x_2x_1^4+x_3^5 + ux_1^2x_2^3
+vx_0x_3K(x_0,...,x_3)=0 \}, \quad u,v \in {\Bbb C}$$ 
of quintic surfaces over $\Spec(\Bbb C[u,v])$ where
$ K$ is a homogenous polynomial of degree 3 with coefficients in 
$\Bbb C[u,v]$. Then there exist elements $\alpha_{u,v}$ in $CH^3(X_{u,v},2)$ 
supported on $Z=X\cap \{x_0x_3=0\}$
such that, for $u,v \in {\Bbb C}$ and $K$ very general, 
these elements are indecomposable modulo the image of $Pic(X_{u,v}) \otimes K_2(\bC)$.
\rm \vskip 6pt

Note that the parameter $v$ is redundant since it is contained in the 
coefficients of the cubic form $K$ already. However in the proof we fix 
$K=x_0^2x_1+x_0x_1x_2+x_0x_2^2+x_2^3+x_0x_1x_3$ and vary $u$ and 
$v$ in a 2-dimensional local parameter space to compute the infinitesimal data.
We suppress $K$ in the notation most of the time.
 
Our examples are deformations of the following quintic hypersurface in
${\Bbb P}^3$:  
$$ 
X:= \{ (x_0:...:x_3) \in {\Bbb P}^3 
\mid F(x)=x_0^5 + x_1x_2^4+x_2x_1^4+x_3^5=0 \}.
$$
The gradient of $F$ is $(5x_0^4:x_2^4+4x_1^3x_2:x_1^4+4x_2^3x_1:5x_3^4)$. It
is nowhere zero, therefore $X$ is smooth. Now we cut out the hyperplane sections
$Z_1:=X \cap \{x_3=0\}$ and $Z_2:=X \cap \{x_0=0\}$. Again $Z_1$ and $Z_2$ are
smooth since the gradients of $Z_1,Z_2$ have no zeroes. The intersection
$$
Z_1 \cap Z_2 = \{ (x_1:x_2) \in {\Bbb P}^1 \mid x_1x_2(x_2^3+x_1^3)=0 \}
$$ 
consists of $5$ distinct points 
$P_1=(0:0:1:0),P_2=(0:1:0:0),P_3=(0:-1:1:0),
P_4=(0:-\zeta:1:0)$ and $P_5=(0:-\zeta^2:1:0),$
where $\zeta^3=1$. The automorphism group
of $X$ contains two copies of ${\Bbb Z}/5{\Bbb Z}$ which are generated by the
diagonal matrices $\sigma_1:=(1:\eta:\eta:\eta)$ and 
$\sigma_2:=(\eta:\eta:\eta:1)$ respectively with $\eta^5=1$. 
$\sigma_i $ operates on $Z_i$ and fixes $Z_{2-i}$ pointwise.  
The fix points of each operation of $\sigma_i$ on $Z_i$ are exactly the
points $P_1,...,P_5$. By the Hurwitz formula the quotient of $Z_i$ by
${\Bbb Z}/5{\Bbb Z}$ 
is a rational curve and the quotient map is ramified in these five
points. Hence on both curves we have the relations
$$ 
5P_1=5P_2=5P_3=5P_4=5P_5 \in CH^1(Z_i),
$$
in particular, on both curves, $5(P_1-P_2)$ and $5(P_3-P_4)$ are rationally
equivalent. On $Z_1$ choose the function $f_1=z$, where $z=x_1/x_2$ is
the pullback of the standard coordinate function on 
$Z_1/\sigma_1 \cong {\Bbb P}^1$. Then $div(f_1)=5(P_1-P_2)$ and $f_1(P_3)=-1,
f_1(P_4)=-\zeta$. Let $g_1:={z+1 \over z+\zeta^2}$. Then
$div(g_1)=5(P_3-P_4)$, $g_1(P_1)=\zeta, g_1(P_2)=1$ and the tame symbol
is 
$$
T(f_1,g_1)=(\zeta^{5},1,-1,-\zeta^{-5},1) \in \oplus_{i=1}^ 5 {\Bbb C}^*.
$$
In a similar way, we can construct $(f_2,g_2)$ on $Z_2$ with the property
that $T(f_2,g_2)=T(f_1,g_1)^{-1}$. 
We see that this element is already quite nice because 
it is non-zero in $CH^2(Z,2)^{ind}$,
but it is 6-torsion. In order to get a non-torsion element, we
have to deform $X$ slightly: look at the two-parameter family
$$
X_{u,v}=\{x_0^5 + x_1x_2^4+x_2x_1^4+x_3^5 + ux_1^2x_2^3
+vx_0x_3K(x_0,...,x_3)=0 \} \quad u,v \in {\Bbb C},
$$ 
where $K$ is a fixed cubic form. 
On all $X_{u,v}$ and their hyperplane sections
$Z_1=\{x_3=0\}$ and $Z_2=\{x_0=0\}$,
we have the same automorphisms $\sigma_1,\sigma_2$ of order $5$, and 
again $Z_1 \cap Z_2$
consists of $5$ fixed points $P_1=(0:0:1:0),P_2=(0:1:0:0),P_3=(0:\alpha:1:0),
P_4=(0:\beta:1:0)$ and $P_5=(0:\gamma:1:0),$ where $\alpha,\beta,\gamma$ are
the $3$ roots of the polynomial $f(z)=z^3+uz+1 \in 
{\Bbb C}[z]$. Define on $C_1$
the rational functions $f_1=z$ and $g_1={z-\alpha \over z-\beta}$. 
Then $f_1(P_3)=\alpha, f_1(P_4)=\beta, f_1(P_5)=\gamma$ and $g_1(P_1)=
{\alpha \over \beta}, g_1(P_2)=g_1(\infty)=1$. 
The tame symbol is given by 
$$
T(f_1,g_1)=(({\beta \over \alpha})^5,1,\alpha^5,\beta^{-5},1) \in 
\oplus_{i=1}^ 5 {\Bbb C}^*.
$$
Now, if $u$ is sufficiently generic (any $v,K$), this tame symbol is not a
torsion element in $\oplus_{i=1}^ 5 {\Bbb C}^*$. We have therefore obtained a 
non-torsion example in $CH^2(Z,2)^{ind}$. 
\vskip 6pt
We now turn to the infinitesimal computation which shows that the cycles
are indecomposable in our sense.
Let $F_{u,v}=x_0^5+x_1x_2^4+x_2x_1^4+x_3^5+ux_1^2x_2^3+vx_0x_3K(x_0,...,x_3)$
be the quintic polynomial, where $K$ is a generic cubic form
with coefficients in $\Lambda=\Bbb C[u,v]$. Let 
$B=\Lambda[x_0,x_1,x_2,x_3]$ and consider the Jacobian rings (cf. \cite{AS1})
$$ \eqalign{&R_F=B/ (\frac{\partial \Fuv}{\partial x_i}\sp (0\leq i\leq 3) ),\cr
& \tRF=B/ ( \frac{\partial \Fuv}{\partial x_i}\sp (i=1,2),\sp
x_j\frac{\partial \Fuv}{\partial x_j}\sp (j=0,3) ). \cr}$$
Then we have isomorphisms (cf. \cite{G2})
$$ H^{2,0}(U/S)\isom \tRF^3=B^3, \sp H^{1,1}(U/S)\isom \tRF^8,$$
$$ H^{2,0}(X/S)\isom R_F^1=B^1,\sp H^{1,1}_{pr}(X/S)\isom R_F^6,$$
and the natural map $H^{p,2-p}(X/S)\to H^{p,2-p}(U/S)$ is given by  
multiplication with $x_0x_3$. The map 
$$\nabla_U:\Omega_S\otimes H^{2,0}(U/S)\to \Omega^2_S\otimes H^{1,1}(U/S)$$
is given by 
$$ Pdu+Qdv \to (P\frac{\partial \Fuv}{\partial v}-Q
\frac{\partial \Fuv}{\partial u})du\wedge dv=
(Px_0x_3K-Qx_1^2x_2^3)du \wedge dv, \sp (P,Q\in B^3). $$
Thus the element lies in $\Ker(\nabla_U)$ if and only if it satisfies 
$$ x_1^2x_2^2(Px_1-Qx_2)\in
\left(\frac{\partial \Fuv}{\partial x_1},
\frac{\partial \Fuv}{\partial x_2},
x_j\frac{\partial \Fuv}{\partial x_j}=5x_j^5\sp (j=0,3)\right).
\leqno{(5-1)}
$$

On the other hand we have the ``residue map'' (5-8)
$$ H^{2,0}(U/S)(= B^3) \to \bigoplus_{i=1}^5 \Lambda;
P\to (Res_{P_i} \frac{\bar{P}}{\bar{F}} \Omega)_{1\leq i\leq 5},$$
where 
$\bar{P}=P\mod <x_0,x_3>\in \Lambda[x_1,x_2]$ and 
$\Omega=x_1dx_2+x_2dx_1=x_2^2 dz\sp (z=\frac{x_1}{x_2})$
is the fundamental form on $\Bbb P^1_{\Lambda}=\Proj\Lambda[x_1,x_2]$. 
The ``residue" of our element is given by 
$$ (d\log \phi_i)_{1\leq i\leq 5}=
\frac{1}{\phi_i}(\frac{\partial \phi_i}{\partial u} du+
\frac{\partial \phi_i}{\partial v} dv)_{1\leq i\leq 5},$$
where
$\phi_1=\beta^5\alpha^{-5},\phi_2=1,\phi_3=
\alpha^5,\phi_4=\beta^{-5},\phi_5=1$.
Thus we want to show the following fact:
there exist no $P,Q\in B^3$ satisfying $(5-1)$ and for $1\leq i\leq 5$
$$ Res_{P_i} \frac{\bar{P}}{\bar{F}} \Omega=
 \frac{1}{\phi_i}\frac{\partial \phi_i}{\partial u} ,\sp
 Res_{P_i} \frac{\bar{Q}}{\bar{F}} \Omega=
 \frac{1}{\phi_i}\frac{\partial \phi_i}{\partial v}.$$

Note that the map $H^{2,0}(U)\to \bigoplus_{i=1}^5 \Lambda$ is 
surjective. Thus our equation
$$ 
Res_{P_i} \frac{\bar{P}}{\bar{F}} \Omega=
 \frac{1}{\phi_i}\frac{\partial \phi_i}{\partial u} ,\sp
 Res_{P_i} \frac{\bar{Q}}{\bar{F}} \Omega=
 \frac{1}{\phi_i}\frac{\partial \phi_i}{\partial v}
 \quad (1\leq i\leq 5)
$$
must have a solution. It is the condition $(5-1)$ 
that destroys such solutions and we show this in the sequel.
Let $a,b,c$ (for simplicity instead of $\alpha,\beta,\gamma$)
be the 3 roots of the equation $f(z)=z^3+uz+1=0$. They satisfy
$abc=-1$, $u=ab+bc+ba$ and $-a-b-c=0$. As above let
$\phi_1=b^5a^{-5},\phi_2=1,\phi_3=a^5,\phi_4=b^{-5},\phi_5=1$.
Using implicit differentiation it is straightforward to prove that
$$
\frac{1}{a}\frac{\partial a}{\partial u}=\frac{-1}{f'(a)}=
\frac{-1}{(a-b)(a-c)}
$$ 
and 
$$
\frac{1}{a}\frac{\partial a}{\partial v}=0,
$$
(since roots do not depend on $v$). 
Therefore we obtain the following set of equations:
$$ \eqalign{
&\frac{1}{\phi_1}\frac{\partial \phi_1}{\partial u}= \frac{-5}{(b-c)(b-a)}+
\frac{5}{(a-b)(a-c)}    \cr
&\frac{1}{\phi_2}\frac{\partial \phi_2}{\partial u}=  0 \cr
&\frac{1}{\phi_3}\frac{\partial \phi_3}{\partial u}=  \frac{-5}{(a-b)(a-c)}\cr
&\frac{1}{\phi_4}\frac{\partial \phi_4}{\partial u}=  \frac{5}{(b-c)(b-a)}\cr
&\frac{1}{\phi_5}\frac{\partial \phi_5}{\partial u}=  0 }
\quad \eqalign{
&\frac{1}{\phi_1}\frac{\partial \phi_1}{\partial v}= 0    \cr
&\frac{1}{\phi_2}\frac{\partial \phi_2}{\partial v}=  0 \cr
&\frac{1}{\phi_3}\frac{\partial \phi_3}{\partial v}=  0 \cr
&\frac{1}{\phi_4}\frac{\partial \phi_4}{\partial v}=  0 \cr
&\frac{1}{\phi_5}\frac{\partial \phi_5}{\partial v}=  0 \cr
}$$
Now we make the convention that the variable $z$
is equal to $\frac{x_1}{x_2}$ and the ``Ansatz''
$$ {\bar P}=a_0+a_1z+a_2z^2+a_3z^3, \sp {\bar Q}=b_0+b_1z+b_2z^2+b_3z^3.
$$
Recall that ${\bar F}=z+uz^2+z^4$ and $\Omega=x_2^2 dz$. 
We then obtain
$$ \eqalign{
&\frac{\bar P}{\bar F}\Omega=\frac{a_0+a_1z+a_2z^2+a_3z^3}{z+uz^2+z^4}dz, \cr
&\frac{\bar Q}{\bar F}\Omega=\frac{b_0+b_1z+b_2z^2+b_3z^3}{z+uz^2+z^4}dz. \cr
}$$
Now compute residues at $P_i$ i.e. $z=0,a,b,c,\infty$:
$$ \eqalign{
Res_{z=0}\frac{\bar P}{\bar F}\Omega &= a_0 \cr
Res_{z=a}\frac{\bar P}{\bar F}\Omega &= 
\frac{a_0+a_1a+a_2a^2+a_3a^3}{a(a-b)(a-c)} \cr
Res_{z=b}\frac{\bar P}{\bar F}\Omega &= 
\frac{a_0+a_1b+a_2b^2+a_3b^3}{b(b-a)(b-c)} \cr
Res_{z=c}\frac{\bar P}{\bar F}\Omega &= 
\frac{a_0+a_1c+a_2c^2+a_3c^3}{c(c-a)(c-b)} \cr
Res_{z=\infty}\frac{\bar P}{\bar F}\Omega &= -a_3 }
\quad \eqalign{
Res_{z=0}\frac{\bar Q}{\bar F}\Omega &= b_0 \cr
Res_{z=a}\frac{\bar Q}{\bar F}\Omega &= 
\frac{b_0+b_1a+b_2a^2+b_3a^3}{a(a-b)(a-c)} \cr
Res_{z=b}\frac{\bar Q}{\bar F}\Omega &= 
\frac{b_0+b_1b+b_2b^2+b_3b^3}{b(b-a)(b-c)} \cr
Res_{z=c}\frac{\bar Q}{\bar F}\Omega &= 
\frac{b_0+b_1c+b_2c^2+b_3c^3}{c(c-a)(c-b)} \cr
Res_{z=\infty}\frac{\bar Q}{\bar F}\Omega &= -b_3 \cr
}$$
This leads to the following 10 equations
$$ \eqalign{
a_0 &=\frac{-5}{(b-c)(b-a)}+ \frac{5}{(a-b)(a-c)}    \cr
\frac{a_0+a_1a+a_2a^2+a_3a^3}{a(a-b)(a-c)}&= 0 \cr
\frac{a_0+a_1b+a_2b^2+a_3b^3}{b(b-a)(b-c)}&= \frac{-5}{(a-b)(a-c)} \cr
\frac{a_0+a_1c+a_2c^2+a_3c^3}{c(c-a)(c-b)}&= \frac{5}{(b-c)(b-a)}\cr
-a_3 &= 0 \cr
b_0 &= 0 \cr
\frac{b_0+b_1a+b_2a^2+b_3a^3}{a(a-b)(a-c)} &=0 \cr
\frac{b_0+b_1b+b_2b^2+b_3b^3}{b(b-a)(b-c)} &=0 \cr
\frac{b_0+b_1c+b_2c^2+b_3c^3}{c(c-a)(c-b)} &=0 \cr
-b_3 &= 0 \cr
}$$
This system, though 8 variables and 10 equations, has always non-zero
solutions for all $a,b,c$, for example in the case
$u=0$ where $a=\exp(i\pi/3),b=-1,c=a^2=\exp(-i\pi/3)$, we get the solution
$$(a_0,...,a_3)=
(- \frac{5}{6} i \sqrt{3} - \frac{5}{2}, 0, - \frac{5}{3}i\sqrt{3}, 0)
=-\frac{5}{3}(a+1,0,a(a+1),0), \sp 
(b_0,...,b_3)=(0,0,0,0).
$$
We therefore need to use condition $(5-1)$ in order to destroy such solutions.
Let us return to homogenous coordinates. Our solutions are
$$
{\bar P}=a_0x_2^3+a_1x_1x_2^2+a_2x_1^2x_2, \sp {\bar Q}=0.
$$
The partial derivatives of $F_{u,v}$ are given by
$$ \eqalign{
\frac{\partial F}{\partial x_0}=& 5x_0^4+x_3K+x_0x_3
\frac{\partial K}{\partial x_0}, \cr
\frac{\partial F}{\partial x_1}=& 
x_2^4+4x_1^3x_2+2ux_1x_2^3+vx_0x_3\frac{\partial K}{\partial x_1}, \cr
\frac{\partial F}{\partial x_2}= &
x_1^4+4x_1x_2^3+3ux_1^2x_2^2 +vx_0x_3\frac{\partial K}{\partial x_2}, \cr
\frac{\partial F}{\partial x_3}= & 5x_3^4+x_0K+x_0x_3\frac{\partial K}
{\partial x_3}.
}
$$
We have to check that the polynomial $R=Px_0x_3K-x_1^2x_2^3Q$
is -- for any two liftings $P={\bar P}+x_0P_0+x_3P_3$ and 
$Q={\bar Q}+x_0Q_0+x_3Q_3$ with quadratic polynomials $P_i,Q_j$ --
not contained in the ideal 
$$ I:=
\left( x_0\frac{\partial F}{\partial x_0},\frac{\partial F}{\partial x_1},
\frac{\partial F}{\partial x_2}, x_3\frac{\partial F}{\partial x_3} \right)
$$
considered as an ideal of ${\Bbb C}[x_0,...,x_3]$.
\vskip 6pt
We compute that 
$$
R(x_0,...,x_3)={\bar P}x_0x_3K+x_0^2x_3KP_0+
x_0x_3^2KP_3-x_0x_1^2x_2^3Q_0-x_1^2x_2^3x_3Q_3,
$$
which is an element of the ideal
$$
J:=\left(x_2^3x_0x_3K,x_1x_2^2x_0x_3K,x_1^2x_2x_0x_3K,
x_0^2x_3K, x_0x_3^2K,x_1^2x_2^3x_0, x_1^2x_2^3x_3 \right).
$$
It is sufficient to show that $I \cap J$ has no 
generators in degree $d \le 8$ for a generic choice of $K$.
It is possible to check this with the computer algebra program 
SINGULAR \cite{S} by using the values $u=1,v=1$ and 
$K=x_0^2x_1+x_0x_1x_2+x_0x_2^2+x_2^3+x_0x_1x_3 $. 
Below is a copy of the corresponding Singular session. All the 
generators of the intersection of the two ideals are computed to 
have degree $\ge 9$.
\vskip 8pt
ring r=0,(w,x,y,z),dp; // here we have $w=x_0,x=x_1,y=x_2,z=x_3$.

poly K=w2x+wxy+wy2+y3+wxz; //sparsepoly(3).

int v = 1;

int u = 1;

poly f1=5w5+wz*K+w2z*diff(K,w);

poly f2=y4+4x3y+2*u*xy3+v*wz*diff(K,x);

poly f3=x4+4xy3+3*u*x2y2+v*zw*diff(K,y);

poly f4=5z5+zw*K+wz2*diff(K,z);

ideal i = f1,f2,f3,f4;//defining ideal of $\tilde R^*$.

ideal j = y3zw*K,xy2zw*K, x2ywz*K, w2z*K, wz2*K, x2y3w, x2y3z;//ideal $J$

ideal I = intersect(i,j);//Intersection $I \cap J$.

dim(I); //the dimension of I (projective dimension +1)!

//Now compute minimal resolution with homogenous entries (lres)!

list T = lres(I,0);

print(betti(T),"betti"); //prints the table of the resolution.

int n;

for (n=ncols(T)); n $>=$1; n=n-1)

{ deg(I[n]), homog(I[n]); }

quit;

\vskip 6pt

In order to finish the proof 
that our cycles are indecomposable, we will need 
that the element we obtained in $H^{1,1}_{pr}(X/S)$ is not
in the Picard group for a general deformation, provided $u,v$ and
$K$ are sufficiently general. This will be true, once the
element is not mapped to zero under the map 
$$
W \otimes H^{1,1}_{pr}(X/S) \longrightarrow H^{0,2}(X/S),
$$
where $W \subset H^1(X,T_X)$ is the 2-dimensional
infinitesimal deformation space of our family parametrized by 
$u,v$ and fixed $K$. Since the element is essentially given by $PK$,
we need to show that $NPK$ is not zero in $R_F^{11}$ for some
quintic polynomial $N \in W \subset R^5_F$, which is generated
by the two elements $x_0x_3K$ and $x_1^2x_2^3$.
We do this in the point $u=0$ with $K=w^2x+wxy+wy^2+y^3+wxz$ 
and $N=x_0x_3K$.
The following little SINGULAR \cite{S} program compares the dimension 
of $R_F^{11}$ with $R_F^{11}/(NPK)$ via Hilbert series and shows 
that $NPK \notin R_F^{11}$: 
\vskip 8pt
ring r=(0,a),(w,x,y,z),dp; minpoly=$a^2-a+1$; 

//here we have $w=x_0,x=x_1,y=x_2,z=x_3$.

poly K= w2x+wxy+wy2+y3+wxz;//sparsepoly(3);

poly P = (a+1)*y3+a*(1+a)*x2y;

poly N = zw*K; // linear combination of x2y3 and zw*K;

poly f=w5+z5+xy4+x4y+zw*K;

ideal j = jacob(f);

ideal i=std(j);

print("Hilbert series of $R_F^*$:");

hilb(i,2); 

print("Hilbert series of $R_F^*/(NPK)$:");

ideal k = jacob(f),N*P*K;

ideal l = std(k);

hilb(l,2);

quit;

\vskip 20pt

\head \S6. An Example in $CH^2(X,1)$. \endhead
\vskip 8pt 

In this section we prove Th.(0-3) which 
provides an example of a class in $CH^2(X,1)$ on a 
smooth projective complex surface $X$ which is indecomposable modulo
the image of $CH^1(X) \otimes {\Bbb C}^*$. 
We state the result here again.

\Pr 6-1. \it On the family 
$$ X_{u}= \{ (x_0:...:x_3) \in {\Bbb P}^3 
\mid F_u(x)=x_0x_1^4+ x_1x_2^4+x_2x_0^4+ x_3^5+ ux_3x_1^4=0 \}
, \quad u \in {\Bbb C}
$$ 
of quintic surfaces, there exist elements $\alpha_u$ in $CH^2(X_{u},1)$ 
such that, for $u$ very general, 
these elements are indecomposable modulo the image of $Pic(X_{u}) \otimes \bC^*$.
\rm \vskip 6pt

In order to construct the examples, we
consider the following quintic hypersurfaces in ${\Bbb P}^3$ from
\cite{SMS}: 
$$ 
X_u:= \{ (x_0:...:x_3) \in {\Bbb P}^3 
\mid F_u(x)=x_0x_1^4+ x_1x_2^4+x_2x_0^4+ x_3^5+ ux_3x_1^4=0 \}.
$$
In \cite{SMS} it was shown that on both curves $Z_1=X_u \cap \{x_3=0\}$
and $Z_2=X_u \cap \{x_0=0\}$ the points $P_1=(0:0:1:0)$ and $P_2=(0:1:0:0)$
satisfy $52(P_1-P_2)=0$ in $CH^1(Z_i)$. This defines an element $\alpha \in 
CH^2(X_u,1)$ for all $u$. It is known by a result of Shioda (cf. \cite{SMS})
that the Picard group of $X_u$ has rank one for almost all $u$.
We use now the method of the previous subsection to deduce that these elements
are indecomposable for very general $u$ modulo $CH^1(X_u) \otimes {\Bbb C}^*$.

We work over a parameter space
$S=\Spec(\Lambda)\subset \Spec(\Bbb C[u])$.
Let $ B=\Lambda[x_0,x_1,x_2,x_3]$
and consider the Jacobian rings
$$ \eqalign{&R_F=B/ (\frac{\partial \Fu}{\partial x_i}\sp (0\leq i\leq 3) )\cr
& \tRF=B/ ( \frac{\partial \Fu}{\partial x_i}\sp (i=1,2),\sp
x_j\frac{\partial \Fu}{\partial x_j}\sp (j=0,3) ). \cr}$$
Then we have isomorphisms (\cite{AS1}, \cite{G2})
$$ H^{2,0}(U/S)\isom \tRF^3=B^3, \sp H^{1,1}(U/S)\isom \tRF^8,$$
$$ H^{2,0}(X/S)\isom R_F^1=B^1,\sp H^{1,1}(X/S)\isom R_F^6,$$
and the natural map $H^{p,2-p}(X/S)\to H^{p,2-p}(U/S)$ is given by 
multiplication with $x_0x_3$.
The map $\nabla_U:H^{2,0}(U/S)\to \Omega^1_S\otimes H^{1,1}(U/S)$
is given by 
$$ Gdu \to G\frac{\partial \Fu}{\partial u} du=
G x_3x_1^4 du , \sp (G\in B^3). $$
Thus the element lies in $\Ker(\nabla_U)$ if and only if it satisfies 
$$Gx_3 x_1^4 \in
\left(\frac{\partial \Fu}{\partial x_1},
\frac{\partial \Fu}{\partial x_2},
x_j\frac{\partial \Fu}{\partial x_j}\sp (j=0,3)\right).
\leqno{(6-1)}
$$
On the other hand we have the ``residue map'' (5-8)
$$ H^{2,0}(U/S)(= B^3) \to \bigoplus_{i=1}^2 \Lambda;
G\to (Res_{P_i} \frac{\bar{G}}{\bar{F}} \Omega)_{1\leq i\leq 2},$$
where 
$\bar{G}=G\mod <x_0,x_3>\in \Lambda[x_1,x_2]$ and 
$\Omega=x_1dx_2+x_2dx_1=x_2^2 dz\sp (z=\frac{x_1}{x_2})$
is the fundamental form on $\Bbb P^1_{\Lambda}=\Proj\Lambda[x_1,x_2]$. 
The ``residue" of our element is given by $\pm 52$ at $P_1$ resp. $P_2$.
Thus we have to show the following fact:
there exist no $G \in B^3$ satisfying $(6-1)$ and for $1\leq i\leq 2$
$$ Res_{P_1} \frac{\bar{G}}{\bar{F}} \Omega= 52 ,\sp
Res_{P_2} \frac{\bar{G}}{\bar{F}} \Omega= -52.
$$ 
It has to be taken also into account 
that $Z_1 \cup Z_2$ is not a normal crossing
divisor, since $P_2$ is a point of multiplicity 4 on the intersection. This gives rise
to 3 more condition which will be discussed below.
Now we make the convention that the variable $z$
is equal to $\frac{x_1}{x_2}$ and the ``Ansatz''
$$ {\bar G}=a_0+a_1z+a_2z^2+a_3z^3. 
$$
Recall that ${\bar F_u}=x_1x_2^4=zx_2^5$ and $\Omega=x_2^2 dz$. 
Therefore we obtain

$$ 
\frac{\bar G}{\bar F}\Omega=\frac{a_0+a_1z+a_2z^2+a_3z^3}{z}dz.
$$
Now compute residues at $P_i$ ($P_2$ is a point of multiplicity 4):
$$ \eqalign{
Res_{z=0}\frac{\bar G}{\bar F}\Omega &= a_0, \cr
Res_{z=\infty}\frac{\bar G}{\bar F}\Omega &= -a_0 }.
$$

This leads to the solution
$$ a_0=52,a_1=a_2=a_3=0. 
$$

We therefore need to use condition $(6-1)$ in order to destroy such solutions.
Let us return to homogenous coordinates. Our solutions are
$$
{\bar G}=a_0x_2^3=52x_2^3.
$$
The partial derivatives of $F_{u}$ are given by
$$ 
\frac{\partial F}{\partial x_0}= x_1^4+4x_0^3x_2,
\frac{\partial F}{\partial x_1}= 4x_0x_1^3+x_2^4+4ux_1^3x_3, 
\frac{\partial F}{\partial x_2}= 4x_1x_2^3+x_0^4, 
\frac{\partial F}{\partial x_3}=5x_3^4+4x_1^4. 
$$
We have to check that the polynomial $R=Gx_1^4x_3$
is -- for any lifting $G=52x_2^3+x_0G_0+x_3G_3$ 
and with quadratic polynomials $G_i$ --
not contained in the ideal 
$$ I:= \left( x_0\frac{\partial F}{\partial x_0},\frac{\partial F}{\partial x_1},
\frac{\partial F}{\partial x_2}, x_3\frac{\partial F}{\partial x_3} \right)
$$
considered as an ideal of ${\Bbb C}[x_0,...,x_3]$.
\vskip 6pt
We compute that 
$R(x_0,...,x_3)=(52x_2^3+x_0G_0+x_3G_3)x_1^4x_3,$
is an element of the ideal
$$
J:=\left( x_1^4x_3  \right).
$$
It is sufficient to show that $I \cap J$ for $u=0$ has only one generator 
$x_0x_1^4x_3$ in degree 6, which does not divide $G$,b and otherwise has no 
generators in degree $d \le 8$.
It is possible to check this with the computer algebra program SINGULAR \cite{S}
by using the value $u=0$. Here is the corresponding Singular session:
\vskip 8pt

ring r=0,(w,x,y,z),dp; //here we have $w=x_0,x=x_1,y=x_2,z=x_3$.

int u = 0;

poly f1=wx4+4*w4*y;

poly f2=4*x3*w+y4+4*u*x3*z;

poly f3=4*x*y3+w4;

poly f4=5*z5+4*x4*z;

ideal i = f1,f2,f3,f4;

ideal j =52*x4y3z, wx4z, x4z2;

ideal I = intersect(i,j);

dim(I); //the dimension of I (projective dimension +1)!

//I; Now compute minimal resolution with homogenous entries (lres)!

list T = lres(I,0);

print(betti(T),"betti");

print(I[1]); print(I[2]);print(I[3]);

quit;

\vskip 20pt

\head \S7. Appendix (by Alberto Collino). \endhead 
\vskip 8pt 

In this section we prove Th.(0-4) which was 
provided to us by Alberto Collino in a letter
from September 19, 1999. We are very grateful to him for letting us reproduce
the contents here. His result shows in particular that indecomposable
cycles in $CH^3(S,2)$ need not be rigid on the surface $S$:

\Pr 7-1. \it On every (very) general quartic $K3$-surface $S$, there exists a
1-dimensional family  of elements $Z_t$ in $CH^3(S,2)$ such that, for $t$ very
general, these elements are indecomposable modulo the image of $Pic(S)
\otimes K_2(\bC)$.\rm

\rm \pr All the cycles which are used originate from the existence
of smooth bielliptic hyperplane sections $C$ of genus $3$ on $S$, which means
that there exists a double cover $C \to E$ onto a smooth elliptic curve $E$
Their existence is guaranteed by the following lemma:

\Lem 7-2. \it On every general quartic surface $S$ there exists a
$1$-dimensional family of bielliptic curves $C_t$ such that the underlying
family $E_t$ of elliptic curves has varying $j$-invariant.
\rm \pr The family $ \Cal A$ of plane quartic curves in $ \Bbb P ^3 $ 
is a projective bundle over the dual projective space $ \Bbb P ^{3 \ast } $. 
Let  $ \Cal B \subset \Cal A$ be the closure of the locus where
the plane quartics are bielliptic, thus $ \Cal B$ is a subvariety of codimension $2$ in 
$\Cal A$.  Let $\Cal K$ be  the linear system of quartic surfaces in $ \Bbb P ^3 $.
Then the incidence family $ \Cal F$ $ \subset \Cal A \times \Cal K $
parametrizes couples $ (Y,S) $ with $  Y \subset S $.  Furthermore, if  $ \Cal Y \subset \Cal F$ 
denotes the restriction of $\Cal F$ over $ \Cal B$, then one has
$ \dim \Cal Y = \dim \Cal K  +1$.
Our aim is to prove that the general quartic surface contains
a 1-dimensional family of bielliptic curves with varying elliptic curve.  
To show this it is enough to check that the tangent space
to a point $ (Y,S) $ in a nonempty fiber of $ \Cal Y \to \Cal K $
is at most of dimension $1$ in general and that there is variation
of the elliptic curve. A smooth bielliptic curve $Y$ has an involution 
acting on it which also acts on the space of canonical 
forms on the curve. The eigenspace with eigenvalue 1 is of dimension
1 while the eigenspace with eigenvalue -1 is of dimension
2.  In the associated coordinates 
the polynomial  of the bielliptic curve is biquadratic in $y$, say it is of the form
$ F(x,y,z) := y^4 +  P(x,z) y^2+ Q(x,z) $.  Now $Y$ is the double cover 
of the elliptic curve $E$ of weighted equation
$  t^2 +  P(x,z) t+ Q(x,z) = 0 $.
Note that $ Y \to E$ is ramified at $y= Q(x,z) = 0$,
while $ E \to \Bbb P^1 $ is ramified over $ P^2 - 4 Q = 0 $.
We remark that a deformation of type 
$ F(x,y,z)+ s G(x,z) =0 $ changes the elliptic curve
if $G(x,z)$ is not in the ideal generated by 
$\partial (P^2 - 4 Q) / \partial x$
and $\partial (P^2 - 4 Q) / \partial z$.
The tangent space $TB$ at $Y$ to the moduli space of bielliptic plane quartic
up to projective equivalence  is isomorphic to the quotient 
of the vector space of polynomials of type
$A(x,z) y^2+ B(x,z) $ modulo the relations
generated by $(\partial F / \partial x) = 0$ and $(\partial F / \partial z) = 0$.
$TB$ is of dimension $4$, if we take
$P := xz $ and $Q := x^ 4 - z^ 4 $ then 
$\{ x^4,x^3z, x^2z^2,xz^3 \}$ is a basis.
The tangent space at Y  to the moduli 
of quartic curves in the plane  
is isomorphic to $R^4$, the degree $4$ part  of the jacobian ring 
of $F$,  with the basis $\{ x^2z y , xz^2 y, x^4,x^3z, x^2z^2,xz^3 \}$.
Consider the quartic surface $S$ with equation 
$$ 
A_2(x,y,z,w)w^2 +A_1(x,y,z)w + F(x,y,z) = 0. 
$$ 
The plane $\pi$ with equation
$w=0$ cuts $S$ along our bielliptic curve $Y$.
Let $w = \epsilon L(x,y,z) $
be the equation of a first order deformation $\pi(L)$  of   $\pi$.
By projecting to $w=0$, the section $S \cap \pi(L)$
gives the curve $Y(L)$ with equation $F + \epsilon LA_1$.
Using the preceeding basis for $R^4$,
one can check that for a given $A_1$ (we have used  
$ A_1 := a x^ 3 + b z^ 3 + y^ 3$ for our computation), the request
that  $Y(L)$ is bielliptic imposes two independent conditions
on $L$ and thus the tangent space at
$Y$ to the family of bielliptic curves on $S$ 
is of dimension $1$ in general.
Furthermore the conditions for the variation of the related elliptic curve 
are also satisfied.
\qed
\vskip 6pt
To continue with the proof of $(7-1)$, let $S$ be a quartic surface in $\bP^3$
such that $S$ admits a smooth bielliptic hyperplane section $C$ of genus $3$
with a double cover $C \to E$ onto a smooth
elliptic curve $E$. Given any element in $CH^2(E,2)$, we can take its
pull-back in $CH^2(C,2)$, which defines an element in $CH^3(S,2)$ by the Gysin
map. In \cite{C} it has been shown that on a
general elliptic curve there are cycles in $CH^2(E,2)$ which
are associated to a pair of Weierstrass points and have a non-trivial
invariant with regard to the regulator map
$$
CH^2(E,2) \to H^1(E,\bC^*).
$$
when $E$ varies in moduli. For convenience let us recall the definition
of those cycles in $CH^2(E,2)$: first choose a non-trivial divisor class
$e$ of degree zero and order $4$ on $E$. Set $\sigma=2e$ and fix a base
point $p \in E$. Then define $q:=p+\sigma$, $a:=p+e$ and $b:=p-e$.
Then there is a unique rational function $f$ on $E$ such that
$div(f)=2q-2p$ and $f(a)=f(b)=1$. Furthermore, by symmetry, there
is a rational function $g$ on $E$ with the property that
$div(g)=2b-2a$ and $g(p)=g(q)=1$. Then the resulting cycle
$\{f,g\} \in CH^2(E,2)$ associated to the graph of
$(f,g); E \to \bP^1 \times \bP^1$ defines a cubical cycle in the
higher Chow group $CH^2(E,2)$.

Now we will do the infinitesimal computations necessary to show proposition
$(7-1)$. Since we use essentially the same methods as in the previous
paragraph, we do not have to introduce any new methods. Let
$S=\{F(x_0,\ldots,x_3)=0 \}$ be the equation of $S$ and $C=\{ x_3=0\}$ the
equation of the hyperplane section $C$. Denote by $U=X \setminus C$ the open
complement. Denote by $T$ a parameter space of pairs $(S,C)$,
where $C \subset S$ is a bielliptic hyperplane section.
Consider, over $T$, the higher Chow group
$CH^{3}(\Cal S,2)$ with support in the divisor $\Cal C$ (the universal
bielliptic section), and inside there the cycle 
$\Cal Z$ which is the pullback of the elliptic cycle from \cite{C}.
The infinitesimal invariant of $\Cal Z$ in our sense
is the lift of the infinitesimal invariant
from the family of elliptic curves,
hence it gives by lemma $(7-2)$ and \cite{C} a non zero element 
$$
\gamma  \in Ker(\Omega^1_{T,0} \otimes H^0(S,\Omega^2_S(\log C))/W_0 
\to \Omega^2_{T,0} \otimes H^1(S,\Omega^1_S(\log C))/W_0). 
$$
Our goal is to prove that any lifting of $\gamma$ into
$H^0(S,\Omega^2_S(\log C))$ is not contained in the kernel of $\nabla_U$.
To this aim we must identify everything 
in terms of Jacobian rings with the previous notations,
$H^{2,0}(U) \cong \tilde R^{1}_F$, and
$H^{1,1}(U) \cong \tilde R^{4}_F$.
After restricting to a $1$-parameter family of surfaces, 
$\Omega^{1}_{T,0}$ has two generators, called $d \sigma$ and $d \tau$, 
dual to the tangent directions 
$\partial \over \partial \sigma$ and $\partial \over \partial \tau$.
$\tilde R^{4}_F$ is the space of tangent directions
to the deformation space of couples $(S,C)$, hence
our vectors $\partial \over \partial \sigma$ and 
$\partial \over \partial \tau$ are associated with polynomials
of degree $4$.  
We now assume that $\partial \over \partial \tau$ corresponds
to a deformation direction where the plane section $C$ is preserved,
and then the polynomial associated to $\partial \over \partial \tau$ 
is of type $\tau := x_3 W$, $W$ a polynomial of degree $3$.
We assume furthermore that $\partial \over \partial \sigma$ corresponds
to a deformation direction where $S$ is preserved. 
Therefore the polynomial associated to $\partial \over \partial \sigma$ 
is of the type $\sigma := L {\partial F \over {\partial x_3}}$, with $L$ a
linear polynomial. The map 
$\nabla_U$ acts on $ d\sigma \otimes P  + d\tau \otimes Q$ by
sending it  to $(d\sigma \wedge d\tau)  ( P \tau  - Q \sigma)$.
Recall that $P$ and $Q$ are polynomials of degree $1$
and that we are working in the jacobian ring $\tilde R^*_F$.

Our element $\gamma$,
the infinitesimal invariant of elliptic origin,
is by its nature of type $ d\sigma \otimes M$,
$M$ of degree $1$ in $x_0, \dots , x_2$.
Indeed in direction $\tau$ the curve is fixed.
So we want to see that, by a wise choice, an element like 
$ d\sigma \otimes M$ is not in the image of the
kernel of $\nabla_U$. Elements which map to $ d\sigma \otimes M$
must be of type $d\sigma \otimes (cx_3  + M) +d\tau \otimes (k x_3)$.
This element is in the kernel of $\nabla_U$
iff in the jacobian ring $\tilde R^*_F$ it is
$(cx_3  + M)x_3 W - (k x_3) L {\partial F \over {\partial x_3}}$ $= 0$,
or, equivalently, iff 
$(cx_3  + M)x_3 W = 0$ in  $\tilde R^5_F$.

Here comes the wise choice, (recall our freedom of choice
of $W$, this amounts to the choice of how the surface changes, keeping
the curve fixed): look at  $(cx_3  + M)x_3$ in  $\tilde R^2_F$, this is clearly not $0$,
hence by duality there is a polynomial $G$ of degree $7$
so that $(cx_3  + M)x_3 G \neq 0 $ in $\tilde R^9_F$.
This implies that there is a polynomial $W$ of degree $3$
so that $W (cx_3  + M)x_3 \neq 0 $ in $\tilde R^5_F$ and
we are done. \qed

\vskip 20pt

\nologo


\enddocument